\documentclass{amsart} 
\usepackage{amssymb,latexsym} 
\usepackage{graphicx, color}
\usepackage{enumerate}

\newcommand{\Row}{\operatorname{Row}}
\newcommand{\Col}{\operatorname{Col}}
\newcommand{\Trace}{\operatorname{Trace}}

\newcommand{\sgn}{\operatorname{sgn}}

\newcommand{\Pictures}{{\rm Pictures}}
\newcommand{\row}{\operatorname{row}}

\newcommand{\QQ}{\mathcal{Q}}

\newcommand{\W}{\mathcal{W}}
\newcommand{\SSS}{\mathcal{S}} 

\newcommand{\Sym}{\mathfrak{S}}

\newcommand{\ZZ}{\mathbb{Z}}
\newcommand{\CC}{\mathbb{C}}

\theoremstyle{plain} 
\newtheorem{theorem}{Theorem}[section] 
\newtheorem{corollary}[theorem]{Corollary}
\newtheorem{lemma}[theorem]{Lemma}
\newtheorem{proposition}[theorem]{Proposition}
\newtheorem{definition}[theorem]{Definition}
\newtheorem{hypotheses}[theorem]{Hypotheses}
\newtheorem{remark}[theorem]{Remark}
\newtheorem{example}[theorem]{Example}
\newtheorem{conjecture}[theorem]{Conjecture}
\newtheorem*{theorem*}{Theorem}

\numberwithin{equation}{section}

\begin{document} 
\title[Skew Schur coincidences] 
{Coincidences among skew Schur functions} 
\subjclass[2000]{Primary 05E05, 20C30} 
\keywords{Symmetric function, skew Schur function, ribbon Schur function, Weyl module}

\thanks{
The first author was supported by NSF grant DMS-0245379.
The second and third authors were supported in part by the
National Sciences and Engineering Research Council of Canada.
The third author was supported in part by 
the Peter Wall Institute for Advanced Studies.}
 \maketitle 

{\parindent 0pt
\emph{author:} Victor Reiner\\
\emph{address:} School of Mathematics\\
University of Minnesota\\
Minneapolis, MN 55455\\
USA\\
\emph{email:} reiner@math.umn.edu
\smallskip\\
\emph{author:} Kristin M. Shaw\\
\emph{address:} Department of Mathematics\\
University of British Columbia\\
Vancouver, BC V6T 1Z2\\
Canada\\ 
\emph{email:} krishaw@math.ubc.ca
\smallskip\\
\emph{author:} Stephanie van Willigenburg\\
\emph{address:} Department of Mathematics\\
University of British Columbia\\
Vancouver, BC V6T 1Z2\\
Canada\\ 
\emph{email:} steph@math.ubc.ca

\medskip

\textbf{All correspondence should be sent to:}\\
\\
Stephanie van Willigenburg\\
Department of Mathematics\\
University of British Columbia\\
Vancouver, BC V6T 1Z2\\
Canada\\
steph@math.ubc.ca}

\pagebreak

\begin{abstract} 
New sufficient conditions and necessary conditions are developed for two skew diagrams to
give rise to the same skew Schur function. 
The sufficient conditions come from a variety of new operations related to {\it ribbons}
(also known as {\it border strips} or {\it rim hooks}). 
The necessary conditions relate to the extent of overlap among the rows or among the
columns of the skew diagram.
\end{abstract}
\section{Introduction \label{intro}}

Symmetric functions play an important role in combinatorics, geometry, and representation theory.
Of particular prominence among the symmetric functions are the family of skew Schur functions $s_{\lambda/\mu}$.
For example, when they were introduced by Schur \cite{Schur} over one hundred years ago they were related to the 
irreducible representations of the symmetric group. Most recently they have been 
connected to branching rules for classical Lie groups \cite{LamPostnikovPylyavskyy, Okounkov}, 
and eigenvalues and singular values of sums of Hermitian and of 
complex  matrices \cite{BergeronBiagioliRosas, FominFultonLiPoon, LamPostnikovPylyavskyy} 
via the study of \emph{inequalities} among products of skew Schur functions. 

With this in mind, a natural avenue to pursue is the \emph{equalities} among products of 
skew Schur functions. As we shall see in Section~\ref{reduction-to-connected}, an equivalent formulation of this question  
is the study of all \emph{binomial syzygies} among skew Schur functions, which is a more tractable 
incarnation of a problem that currently seems out of reach: find \emph{all} syzygies among skew Schur functions.
Famous non-binomial syzygies include various formulations of the Littlewood-Richardson rule and 
Equation \eqref{most-basic-syzygy} below,
which give some indication of the complexity that any solution would involve.

The study of equalities among skew Schur functions can also be regarded as part of the ``calculus of shapes''. 
For an arbitrary subset $D$ of $\ZZ^2$, there are two polynomial
representations $\SSS^D$ and $\W^D$ of $GL_N(\CC)$ known as a \emph{Schur module} and \emph{Weyl module} respectively,
obtained by row-symmetrizing and column-antisymmetrizing tensors whose tensor positions are
indexed by the cells of $D$.  These representations are determined up to isomorphism by their
character, namely the symmetric function $s_D(x_1,\ldots,x_N)$, which tells us the trace of
any element $g$ in $GL_N(\CC)$ acting on $\SSS^D$ and $\W^D$ as a function of the eigenvalues $x_1,\ldots,x_N$ 
of $g$.  When $D=\lambda/\mu$ is a skew diagram, this symmetric function is the skew Schur function
$s_{\lambda/\mu}(x_1,\ldots,x_N)$.  Therefore, the question of when two skew Schur or Weyl modules are 
equivalent, working over $\CC$, is precisely the question of equalities among skew Schur functions.

As a consequence of this, the aim of this paper is to study the equivalence relation on skew diagrams $D_1, D_2$ defined by $D_1\sim D_2$ if and only if $s_{D_1}=s_{D_2},$ and in particular to use known \emph{skew-equivalences} to generate new ones. Our motivation for this approach is \cite{BilleraThomasvanWilligenburg} where Billera, Thomas and the third author studied when two elements of the  subclass of skew diagrams known as \emph{ribbons} or \emph{border strips} or \emph{rim hooks} were skew-equivalent. They discovered that if ribbons $\alpha ,\beta, \gamma ,\delta$ satisfied $\alpha \sim \beta$ and $\gamma \sim \delta$ then the \emph{composition} of ribbons $\alpha \circ \beta$ and $\gamma \circ \delta$ satisfied $\alpha \circ \beta \sim\gamma \circ \delta$.

The paper is structured as follows. In Section \ref{diagrams} we review notation concerning partitions, compositions and skew diagrams. 
Section \ref{Lambda} recalls the ring of symmetric functions 
and Section \ref{Schur-function-definitions} covers various definitions and basic properties of skew Schur functions.
Section \ref{Lrrule} is our final review section and gives a version of the Littlewood-Richardson rule.

In Section \ref{reduction-to-connected} we reduce the question of skew-equivalence to the case of 
connected skew diagrams. 
Sections \ref{sufficiencies} and \ref{necessities} then build upon this to develop necessary and sufficient conditions for skew-equivalence. Specifically, in Section \ref{sufficiencies}, for ribbons $\alpha , \beta$ and a skew diagram $D$ we define compositions $\alpha \circ D$ and $D\circ \beta$ that naturally generalise the composition of ribbons, $\circ$, defined in \cite{BilleraThomasvanWilligenburg}  and prove

\begin{theorem*}(Theorem \ref{composition-equivalences})
If one has ribbons $\alpha, \alpha'$ and skew diagrams $D, D'$ satisfying
$\alpha \sim \alpha'$ and $D \sim D'$, then
\begin{enumerate}
\item[(i)]$\alpha \circ D \sim     \alpha' \circ D,$ 
\item[(ii)] $D \circ \alpha \sim D' \circ \alpha,$ 
\item[(iii)] $D \circ \alpha   \sim    D \circ \alpha',$ and
\item[(iv)]$\alpha \circ D \sim    \alpha \circ D^\ast,$
\end{enumerate}
where $D^\ast$ is $D$ rotated by 180 degrees.
\end{theorem*}

For certain ribbons $\omega$ we also construct an analogous operation to $\circ$ called \emph{amalgamated composition}, $\circ _\omega$, and prove

\begin{theorem*}(Theorem \ref{RibbComp})
If $\alpha, \alpha'$ are ribbons with $\alpha \sim \alpha'$, and
 $D, \omega$ satisfy Hypotheses~\ref{wseparation}, then one has the following skew-equivalences:
$$
\alpha' \circ_{\omega} D
\,\, \sim \,\,
\alpha \circ_{\omega} D 
\,\, \sim \,\,
\alpha \circ_{\omega^\ast} D^\ast ,$$where $D^\ast$ is $D$ rotated by 180 degrees.
\end{theorem*}

Additionally, Section \ref{ribbon-staircases} yields a construction 
that produces skew diagrams that are skew-equivalent to their conjugate.

Meanwhile, Section \ref{necessities} discusses two necessary conditions for skew-equivalence. 
One comes from the Frobenius rank of a skew diagram studied in \cite{ChenYang, Stanley2, YanYangZhou}.
The other is new, and relates to the sizes of the rows and the columns of a skew diagram, 
and the sizes of their overlaps. Finally, Section \ref{Conclusion} suggests further avenues to pursue.

\tableofcontents

\subsection{Acknowledgements} 
The authors would like to thank Louis Billera, Peter McNamara, Richard Stanley and John Stembridge for helpful discussions,
Christopher Ryan and John Stembridge for aid with data generation, and all
those who suggested the question. They would also like to thank the referee for insightful and constructive comments that improved the exposition.

%\subsection{Notation table (to be removed later)}
%
%$$
%\begin{array}{ll}
%\alpha,\beta,\gamma,\rho&\mbox{composition/ribbon}\\
%\lambda,\mu,\nu&\mbox{partition/diagram}\\
%\varepsilon ^k_m()&\mbox{ribbon staircase}\\
%\ast&\mbox{antipodal rot/reversing a composition}\\
%m,n&\mbox{composition size}\\
%\cdot ,  \odot ,\circ&\mbox{concat, near concat, comp of comp}\\
%\star&\mbox{mirror of nesting}\\
%l(),|\ |&\mbox{length, size} \text{  ({\bf NB}:  Vic chose to use }\ell\text{ instead of $l$ )}\\
%t&\mbox{conjugate}\\
%\delta&\mbox{staircase}\\
%\cap _m, \cup _m&\mbox{$m$-intersection, -union}\\
%\sqcap&\mbox{overlap set}\\
%\theta \# \theta&\mbox{H-G ribbons}\\
%s&\mbox{Schur function (i.e. not $S$)}\\
%\omega&\mbox{involution}\\
%\Lambda&\mbox{symm fn}\\
%\mathcal{N}&\mbox{nesting}\\
%\mathcal{R,O}&\mbox{outside decomp, coarsening}\\
%\kappa()&\mbox{coarsening factor}\\
%\sim&\mbox{Schur equivalent}\\
%D^+&\mbox{add a column}\\
%r^{(k)},\rho ^{(k)}&\mbox{row overlap/partition}\\
%c^{(k)},\gamma ^{(k)}&\mbox{col overlap/partition}\\
%t_k(),b_k()&\mbox{subribbons}
%\end{array}
%$$

\section{Diagrams \label{diagrams}}

  In this section, we review  partitions,
compositions, Ferrers diagrams, skew diagrams and ribbons.  The interested reader may wish to consult \cite{Macdonald, Sagan, Stanley1} for further details. 

A \emph{partition} $\lambda$ of a positive integer $n$, denoted $\lambda \vdash n$,  is a sequence 
$(\lambda_1,\lambda_2,\ldots,\lambda_\ell)$ of
positive integers $\lambda_i$ such that
$$
\lambda_1 \geq \cdots \geq \lambda_\ell > 0
$$
and $\sum_{i=1}^\ell \lambda_i=n$.  We call $n$ the \emph{weight} or \emph{size}
of $\lambda$, and denote it
$|\lambda|:=n$.  Each $\lambda _i$ is called a
\emph{part} of $\lambda$, and the number of parts $\ell$ is called the
\emph{length} $\ell(\lambda):=\ell$. The unique partition of 0 is denoted by $\varnothing$.

The (\emph{Ferrers} or \emph{Young}) \emph{diagram} of $\lambda$ consists of \emph{boxes} or \emph{cells} such that there are $\lambda_i$ cells  in each
row $i$, so the top row has $\lambda_1$ cells,
the second-from-top row has $\lambda_2$ cells, etc. In addition, the rows of cells are all left-justified. We abuse notation and also denote the Ferrers
diagram of $\lambda$ by $\lambda$.

Two partial orders on partitions that arise frequently are 
\begin{enumerate}
\item[$\bullet$] 
the \emph{inclusion order}: $\mu \subseteq \lambda$ if $\mu_i \leq \lambda_i$ for all $i$,
\item[$\bullet$]
the \emph{dominance} (or \emph{majorization}) \emph{order} on partitions $\lambda, \mu$ having the \emph{same
weight}:  $\mu \leq_{dom} \lambda$ if 
$$
\mu_1 + \mu_2 + \cdots + \mu_i \leq \lambda_1 + \lambda_2 + \cdots + \lambda_i
$$
for $i=1,2,\ldots,\min(\ell(\mu),\ell(\lambda))$.
\end{enumerate}

Given two partitions $\lambda, \mu$ such that $\mu \subseteq \lambda$ the \emph{skew (Ferrers) diagram} $D=\lambda/\mu$
is obtained from the Ferrers diagram of $\lambda$ by removing the cells in the subdiagram
of $\mu$ from the top left corner. 
For example, the following is a skew diagram whose cells are indicated by $\times$:
$$
\lambda/\mu=(5,4,3,3)/(3,1) =
\begin{matrix}
      &      &      &\times&\times \\
      &\times&\times&\times&       \\
\times&\times&\times&      &       \\
\times&\times&\times&      &
\end{matrix}.
$$

Cells in skew diagrams will be referred to by their row and column indices 
$(i,j)$, where $i \leq \ell(\lambda)$ and $j \leq \lambda_i$.  The \emph{content} or
\emph{diagonal} of the cell is the integer $c(i,j)=j-i$.

Given two skew diagrams $D_1, D_2$, a \emph{disjoint union} $D_1 \oplus D_2$
of them is obtained by placing $D_2$ strictly to the north and east of $D_1$ in such a way that
$D_1, D_2$ occupy none of the same rows or columns.  For example, if $D_1 = (2,2), D_2 = (3,2)/(1)$ then a possible disjoint union is 
\begin{equation}
\label{disjoint-sum-example}
D_1 \oplus D_2 = 
\begin{matrix}
      &      &      &\times&\times \\
      &      &\times&\times&       \\
\times&\times&      &      &       \\
\times&\times&      &      &
\end{matrix}.
\end{equation}

We say that a skew diagram $D$ is
is \emph{connected} if it cannot be written as $D=D_1 \oplus D_2$ for two proper
subdiagrams $D_1,D_2$.  A connected skew diagram $D$ is
called a \emph{ribbon} or \emph{border strip} or \emph{rim hook}
if it does not contain a subdiagram isomorphic to that of the partition $\lambda=(2,2).$
For example,
\begin{equation}
\label{example-of-ribbon}
\lambda/\mu=(5,4,3,1)/(3,2) =
\begin{matrix}
      &      &      &\times&\times \\
      &      &\times&\times&       \\
\times&\times&\times&      &       \\
\times&      &      &      &
\end{matrix}
\end{equation}
is a ribbon.  Two
skew diagrams $D, \tilde{D}$ will be considered equivalent 
as subsets of the plane
%\footnote{This notion of equivalence 
%is easily seen to be consistent with the notion of Schur equivalence defined in
%the introduction:  two skew diagrams that are equivalent as subsets of the plane
%will have $D_1 \sim D_2$, that is, $s_{D_1} = s_{D_2}$.}
if one can be obtained from the other by vertical or horizontal translations, or by the removal or addition of empty rows or columns. As a consequence of this, given two diagrams $D_1, D_2$ we can now say   their {disjoint union} $D_1 \oplus D_2$
 is obtained by placing $D_2$ immediately to the north and east of $D_1$ in such a way that
$D_1, D_2$ occupy none of the same rows or columns, as illustrated by Equation~\eqref{disjoint-sum-example}.

 A \emph{composition} $\alpha$ of a positive integer $n$, denoted $\alpha \vDash n$,  is an ordered sequence 
$(\alpha_1,\alpha_2,\ldots,\alpha_\ell)$ of
positive integers $\alpha_i$ such that $\sum_{i=1}^\ell \alpha_i=n$.  
As with partitions, we call $n$  the \emph{weight} or \emph{size} of $\alpha$, and denote it by
$|\alpha|:=n$.  
Again, the number $\ell$ is called the \emph{length} $\ell(\alpha):=\ell$.

 We end with two bijections regarding compositions.  For a positive integer
$n$, let $[n]:=\{1,2,\ldots,n\}$.  For the first bijection consider the map sending a composition 
$\alpha=(\alpha_1,\ldots,\alpha_\ell)$ to the set of partial sums
$\{\alpha_1,\alpha_1+\alpha_2,\ldots,\alpha_1+\alpha_2+\cdots+\alpha_{\ell-1} \}$,
which gives a bijection between compositions of $n$ and the collection 
$2^{[n-1]}$ of all subsets of $[n-1]$.  For the second bijection consider the map sending $\alpha$ to the unique
ribbon having $\alpha_i$ cells in the $i^{th}$ row from the 
\emph{bottom}, which
gives a bijection between compositions of $n$ and ribbons of size $n$. 
Note that labelling the rows of a composition from bottom to top
is slightly inconsistent with the labelling of rows of Ferrers
diagrams from top to bottom in English notation, but it is in keeping with the seminal work \cite{Gelfandetal}.
Due to this bijection, we will often refer to ribbons by their composition of
row sizes.  To illustrate these bijections, observe that
the composition $\alpha=(1,3,2,2)$ of $n=8$ corresponds to the
subset $\{1,4,6\}$ of $[n-1]=[7]$, and to the ribbon depicted in \eqref{example-of-ribbon}.

\subsection{Symmetries of diagrams}

We will have occasion to use several symmetries of partitions and skew diagrams
and review two of them here.

Given a partition $\lambda$, its \emph{conjugate} or \emph{transpose} partition
$\lambda^t$ is the partition whose Ferrers diagram is obtained from that of $\lambda$
by reflecting across the northwest-to-southeast diagonal.  Equivalently, the
parts of $\lambda^t$ are the column sizes
of the Ferrers diagram of $\lambda$ read from left to right.  This extends to skew diagrams in a natural
way:  if $D=\lambda/\mu$ then $D^t:=\lambda^t/\mu^t$.  

Given a skew diagram $D$, one can
form its \emph{antipodal rotation} $D^*$ by rotating it $180$ degrees in the plane.
Note that for a ribbon $\alpha=(\alpha_1,\ldots,\alpha_\ell)$, the antipodal rotation
of its skew diagram corresponds to the \emph{reverse} composition  
$\alpha^*=(\alpha_\ell,\ldots,\alpha_1)$.

\subsection{Operations on ribbons and diagrams}
\label{ribbon-operations}
  This subsection reviews some standard operations on ribbons.  It also discusses
a composition operation $\alpha \circ \beta$ on ribbons $\alpha, \beta$
that was introduced in \cite{BilleraThomasvanWilligenburg}, and its generalization
to operations $\alpha \circ D$ and $D \circ \beta$ for skew diagrams $D$.

  Given two skew diagrams $D_1, D_2$, aside from their disjoint sum $D_1 \oplus D_2$,
there are two closely related important operations called their \emph{concatentation} $D_1 \cdot D_2$
and their \emph{near-concatenation} $D_1 \odot D_2$.  The concatentation $D_1 \cdot D_2$ (resp.
near concatentation $D_1 \odot D_2$) is
obtained from the disjoint sum $D_1 \oplus D_2$ by moving all cells of $D_2$
one column west (resp. one row south), so that the same column (resp. row) is occupied by the
rightmost column (resp. topmost row) of $D_1$ and the leftmost column (resp. bottommost row)
of $D_2$. For example, if 
$$
\begin{aligned}
D_1 &= (2,2)\\
D_2 &= (3,2)/(1)
\end{aligned}
$$
then $D_1 \oplus D_2$ was given in Equation
\eqref{disjoint-sum-example}, while
$$
D_1 \cdot D_2 = 
\begin{matrix}
      &      &2     &2 \\
      &2     &2     &       \\
1     &1     &      &       \\
1     &1     &      &
\end{matrix}
\quad
D_1 \odot D_2 =
\begin{matrix}
      &      &      &2     &2 \\
1     &1     &2     &2     &       \\
1     &1     &      &      &       
\end{matrix}.
$$
Observe we have used the numbers $1$ and $2$ to distinguish the cells in $D_1$ from the cells in $D_2$.  
The reason for the names ``concatentation'' and ``near-concatentation'' becomes
clearer when we restrict to ribbons.  Here if
$$
\begin{aligned}
\alpha&=(\alpha_1,\ldots,\alpha_\ell)\\
\beta&=(\beta_1,\ldots,\beta_m),
\end{aligned}
$$
then 
$$
\begin{aligned}
\alpha \cdot \beta &= (\alpha_1,\ldots,\alpha_\ell,\beta_1,\ldots,\beta_m) \\
\alpha \odot \beta &= (\alpha_1,\ldots,\alpha_{\ell-1},\alpha_\ell+\beta_1,\beta_2,\ldots,\beta_m),
\end{aligned}
$$
which are the definitions for concatenation and near concatenation given in \cite{Gelfandetal}.

 Note that the operations $\cdot$ and $\odot$
are each associative, and associate with each other:
\begin{equation}
\label{associativity}
\begin{aligned}
( D_1 \cdot D_2 ) \cdot D_3 &= D_1 \cdot ( D_2  \cdot D_3 ) \\
( D_1 \odot D_2 ) \odot D_3 &= D_1 \odot ( D_2  \odot D_3 ) \\
( D_1 \odot D_2 ) \cdot D_3 &= D_1 \odot ( D_2  \cdot D_3 ) \\
( D_1 \cdot D_2 ) \odot D_3 &= D_1 \cdot ( D_2  \odot D_3 )
\end{aligned}.
\end{equation}
Consequently a string of operations $D_1 \star_1 D_2 \star_2 \cdots \star_{k-1} D_k$
in which each $\star_i$ is either $\cdot$ or $\odot$ is well-defined
without any parenthesization.   Also note that ribbons are exactly the
skew diagrams that can be written uniquely as a string of the form
\begin{equation}
\label{alpha-ribbon-string}
\alpha = \square \star_1 \square \star_2 \cdots \star_{k-1} \square
\end{equation}
where $\square$ is the diagram with exactly one cell.

Given a composition $\alpha$ and a skew diagram $D$, define $\alpha \circ D$ to be
the result of replacing each cell $\square$ by $D$ in the expression \eqref{alpha-ribbon-string}
for $\alpha$:
$$
\alpha \circ D := D \star_1 D \star_2 \cdots \star_{k-1} D.
$$
For example, if 
$$
\alpha=(2,3,1) = 
\begin{matrix}
      &      &      &\times\\
      &\times&\times&\times\\
\times&\times&      &
\end{matrix}
\quad\text{ and }\quad
D=
\begin{matrix}
\times&\times&\\
\times&\times&
\end{matrix}
$$
then 
$$
\begin{aligned}
\alpha &= \square\odot\square \cdot 
             \square\odot\square\odot\square \cdot 
                   \square \\
\alpha \circ D &=
           D\odot D \cdot 
             D\odot D\odot D \cdot 
                   D \\
&=
\begin{matrix}
 & & & & & & & &6&6\\
 & & & & & & & &6&6\\
 & & & & & & &5&5& \\
 & & & & &4&4&5&5& \\
 & & &3&3&4&4& & & \\
 & & &3&3& & & & & \\
 & &2&2& & & & & & \\
1&1&2&2& & & & & & \\
1&1& & & & & & & &
\end{matrix}
\end{aligned}
$$
where we have used numbers to distinguish between copies of $D$.

  It is easily seen that when $D=\beta$ is a ribbon, then $\alpha \circ \beta$ is also a
ribbon, and  agrees with the definition in \cite{BilleraThomasvanWilligenburg}.

  Similarly, given a skew diagram $D$ and a ribbon $\beta$,
we can also define $D \circ \beta$  as follows.  Create a copy $\beta^{(i)}$ of
the ribbon $\beta$ for each of the cells of $D$, numbered $i=1,2,\ldots,n$ arbitrarily.
Then assemble the diagrams $\beta^{(i)}$ into a disjoint decomposition of
$D \circ \beta$ by translating them in the plane, in such a way that $\beta^{(i)} \sqcup \beta^{(j)}$ forms
a copy of
\begin{equation}\label{D-circ-beta-defn}
\begin{cases}
\beta^{(i)} \odot \beta^{(j)} &\text{ if }i\text{ is just left of }j\text{ in some row of }D, \\
\beta^{(i)} \cdot \beta^{(j)} &\text{ if }i\text{ is just below }j\text{ in some column of }D. \\
\end{cases}\end{equation}

For example, if
$$
D=
\begin{matrix}
 &1&2\\
3&4&5
\end{matrix},
\qquad
\beta =
\begin{matrix}
      &\times&\times&\times\\
     \times&\times& &
\end{matrix}
$$
then $D \circ \beta$ is the skew diagram
$$
\begin{array}{*{15}c}
  & & & & & & & & & & & &2&2&2\\
  & & & & & & & &1&1&1&2&2& & \\
  & & & & & & &1&1&5&5&5& & & \\
  & & & & &4&4&4&5&5& & & & & \\
  &3&3&3&4&4& & & & & & & & & \\
 3&3& & & & & & & & & & & & &
\end{array}
$$
where we have used numbers to distinguish between copies of $\beta$. One must check that the local constraints defining $D \circ \beta$ given 
in \eqref{D-circ-beta-defn} are indeed simultaneously satisfiable 
globally, and hence that $D \circ \beta$ is well-defined.  For this it 
suffices to check the case $D=\lambda = (2,2)$, which we leave to the reader as an easy
exercise.
Again it is clear that when $D=\alpha$ is a ribbon, then $\alpha\circ\beta$ is another ribbon agreeing with that in \cite{BilleraThomasvanWilligenburg}.
The following distributivity properties should also be clear.

\begin{proposition}
\label{circ-distributivities}
For skew diagrams $D, D_1, D_2$ and ribbons $\alpha$ and $\beta$ the operation $\circ$ distributes over  $\cdot$ and $\odot$, that is
$$
\begin{aligned}
(\alpha \cdot \beta) \circ D &= (\alpha \circ D) \cdot (\beta \circ D) \\
(\alpha \odot \beta) \circ D &= (\alpha \circ D) \odot (\beta \circ D) \\
\end{aligned}
$$
and
%\footnote{{\bf Steph:}  Would you check me on the 3rd, 4th
%distributivities in Proposition~\ref{circ-distributivities}?
%I don't think we mentioned them before, because we didn't need them (and still don't).  But
%I think they're true, and it would be sinful not to include them, if true.}
$$
\begin{aligned}
(D_1 \cdot D_2) \circ \beta &= (D_1 \circ \beta) \cdot (D_2 \circ \beta) \\
(D_1 \odot D_2) \circ \beta &= (D_1 \circ \beta) \odot (D_2 \circ \beta) \\
\end{aligned}.
$$
\end{proposition}
\noindent

\begin{remark} \rm \ \\
Observe that $D_1\circ D_2$ has not been defined for both $D_1$ and $D_2$ being non-ribbons,
as certain difficulties arise.  We invite the reader to investigate this already in
the case where $D_1, D_2$ are both equal to the smallest non-ribbon, namely
the $2 \times 2$ rectangular Ferrers diagram $\lambda=(2,2)$, in order to appreciate these difficulties;
 see also Remark~\ref{smallest-nonribbon-circ-example} below.
\end{remark}

\section{The ring of symmetric functions \label{Lambda}}

We now recall the ring of symmetric functions $\Lambda$,
and some of its polynomial generators and bases. Further details can be found in the excellent texts \cite{Macdonald, Sagan, Stanley1}.

  The ring $\Lambda$ is the subalgebra of the formal power series
$\ZZ[[x_1,x_2,\ldots]]$ in countably many variables, consisting of those
series $f$ that are of bounded degree in the $x_i$, and invariant under all
permutations of the variables.  If $\Lambda^n$ denotes the
symmetric functions that are homogeneous of degree $n$, then
we have an abelian group direct sum decomposition $\Lambda = \bigoplus_{n \geq 0}{\Lambda^n}$.
There is a natural $\ZZ$-basis for $\Lambda^n$ given by
the \emph{monomial symmetric functions} $\{ m_\lambda \}_{\lambda \vdash n}$,
where $m_\lambda$ is the formal sum of all monomials that can be permuted
to $x^{\lambda}:=x_1^{\lambda_1} \cdots x_\ell^{\lambda_\ell}$.

  The fundamental theorem of symmetric functions states that
$\Lambda$ is a polynomial algebra in the \emph{elementary symmetric functions}
$$
\Lambda = \ZZ[e_1,e_2,\ldots]
$$
where 
$$
e_r:= \sum_{ 1 \leq i_1 < i_2 < \cdots < i_r } x_{i_1} x_{i_2} \cdots x_{i_r}.
$$
It transpires that it is also a polynomial algebra in the \emph{complete
homogeneous symmetric functions} 
$$
h_r:= \sum_{ 1 \leq i_1 \leq i_2 \leq \cdots \leq i_r } x_{i_1} x_{i_2} \cdots x_{i_r},
$$
and the map $\omega: \Lambda \rightarrow \Lambda$ mapping $e_r \longmapsto h_r$
is an involution.  To obtain $\ZZ$-bases for $\Lambda$, define
for partitions $\lambda=(\lambda_1,\ldots,\lambda_\ell)$
$$
\begin{aligned}
e_\lambda & :=e_{\lambda_1} \cdots e_{\lambda_\ell} \\
h_\lambda & :=h_{\lambda_1} \cdots h_{\lambda_\ell}
\end{aligned}.
$$
From here a consequence of the fundamental theorem is that 
$\Lambda^n$ has as a $\ZZ$-basis either
$\{e_\lambda\}_{\lambda \vdash n}$ or $\{h_\lambda\}_{\lambda \vdash n}$.

\section{Schur and skew Schur functions\label{Schur-function-definitions}}

  This section reviews some definitions of  Schur functions $\{ s_\lambda\} _{\lambda \vdash n , n\geq 0}$ and skew Schur functions that will be useful.

 \subsection{Tableaux}

  One way to define the (skew) Schur function $s_D$ for a (skew) diagram $D$
involves tableaux.  A \emph{column-strict} (or \emph{semistandard}) \emph{tableau} of \emph{shape} $D$
is a \emph{filling} $T: D \rightarrow \{1,2,\ldots\}$ of  the cells of
$D$ with positive integers such that the numbers 
\begin{enumerate}
\item[(i)]
weakly increase left-to-right in each row,
\item[(ii)]
strictly increase top-to-bottom down each column.
\end{enumerate}
The \emph{(skew) Schur function} $s_D$ is then
\begin{equation}
\label{tableau-Schur-definition}
s_D : = \sum_{T} x^T
\end{equation}
where the sum ranges over all column-strict tableaux of shape $D$, and
$$
x^T:= \prod_{(i,j) \in D} x_{T{(i,j)}}.
$$If $D$ is a ribbon we call $s_D$ a \emph{ribbon Schur function}. That (skew) Schur functions are symmetric follows from the definition\begin{equation}
\label{Kostka-numbers}
\begin{aligned}
s_D &= \sum_{\mu} K_{D,\mu} m_\mu  .\\
%\langle s_D , h_\mu \rangle &= K_{D,\mu}
\end{aligned}
\end{equation}
Here $K_{D,\mu}$ is the \emph{Kostka number}, which is
number of column-strict tableaux of shape $D$
and \emph{content} $\mu$, that is, having $\mu_i$ occurrences of $i$ for each $i$.
From the definition \eqref{tableau-Schur-definition}, 
one of the most basic syzygies \cite[Chapter 1.5, Example 21 part (a)] {Macdonald}
among skew Schur functions follows immediately.

\begin{proposition}
\label{most-basic-syzygy} If $D_1$ and $D_2$ are skew diagrams then
$$
s_{D_1} s_{D_2} = s_{D_1 \cdot D_2} + s_{D_1 \odot D_2}.
$$
\end{proposition}
\begin{proof}
Given a pair $(T_1,T_2)$ of column-strict tableaux of shapes $(D_1,D_2)$,
let $a_1$ be the northeasternmost entry of $T_1$ and $a_2$ the
the southwesternmost entry of $T_2$.  Then either
\begin{enumerate}
\item[$\bullet$]
$a_1 > a_2$, and hence $(T_1,T_2)$ concatenate to make a column-strict
tableaux of shape $D_1 \cdot D_2$, or
\item[$\bullet$]
$a_2 \geq a_1$, and hence $(T_1,T_2)$ near-concatenate to make a column-strict
tableaux of shape $D_1 \odot D_2$.
\end{enumerate}
\end{proof}

 \subsection{The Jacobi-Trudi determinant and the infinite Toeplitz matrix}
\label{Toeplitz-section}

Skew Schur functions turn out to be the nonzero minor subdeterminants in
certain Toeplitz matrices.  Consider the sequence ${\bf h}:=(h_0(=1),h_1,h_2,\ldots)$
and its \emph{Toeplitz matrix}, the infinite matrix 
$$
T:=(t_{ij})_{i,j \geq 0} := (h_{j-i})_{i,j \geq 0}
$$
with the convention that $h_r=0$ for $r < 0$.
The \emph{Jacobi-Trudi determinant} formula for the skew Schur function $s_{\lambda /\mu}$ asserts that 
\begin{equation}
\label{Jacobi-Trudi-formula}
s_{\lambda/\mu} = \det(h_{\lambda_i - \mu_j -i + j})_{i,j=1}^{\ell(\lambda)}.
\end{equation}

This can be reinterpreted as follows:
the square submatrix of the Toeplitz matrix $T$
having row indices  $i_1 < \ldots < i_m$ and column indices $j_1 < \ldots < j_m$
has determinant equal to the skew Schur function $s_D$ for
$D=\lambda/\mu$ where for $r=1,2,\ldots,m$ 
\begin{equation}
\label{Toeplitz-minor-to-partitions}
\begin{aligned}
\lambda_r  &:= j_m - i_r - m + r \\
\mu_r      &:= j_m - j_r - m + r.
\end{aligned}
\end{equation}
In particular, if for some $r$ one has $\lambda_r < \mu_r$, then this determinant will be zero.

  We remark here that transposing a skew diagram $D$ 
corresponds to the involution $\omega$ on $\Lambda$ that exchanges $e_r$ and $h_r$ for all $r$, that is
\begin{equation}
\label{omega-is-transpose}
\omega (s_D) = s_{D^t}.
\end{equation}
As a consequence, there is a \emph{dual Jacobi-Trudi determinant} that
is obtained by applying $\omega$ to \eqref{Jacobi-Trudi-formula},
which expresses $s_D$ as a polynomial in the elementary symmetric functions $e_r$.

 \subsection{The Hamel-Goulden determinant}

  One can view the Jacobi-Trudi determinant (or its dual) as expressing a skew Schur function 
in terms of skew Schur functions of particular shapes, namely shapes consisting of a single
row (resp. a single column),
since by the definition  \eqref{tableau-Schur-definition} $h_r=s_r$ (resp. $e_r=s_{1^r}$).  
There are other such determinantal formulae for Schur and skew Schur functions such as 
the \emph{Giambelli determinant} involving hook shapes, 
the \emph{Lascoux-Pragacz determinant} involving ribbons \cite{LascouxPragacz}, 
and most generally the \emph{Hamel-Goulden determinant} \cite{HamelGoulden}.  We review this last 
determinant here, using the reformulation
involving the notion of a \emph{cutting strip} due to Chen, Yan and Yang \cite{ChenYanYang}.

  Given a skew diagram $D$, an  \emph{outside (border strip) decomposition} is an ordered decomposition
$\Pi=(\theta_1,\ldots,\theta_m)$ of $D$, where each $\theta_k$ is a ribbon whose southwesternmost
(resp. northeasternmost) cell lies either on the left or bottom (resp. right or top) 
perimeter of $D$.  Having fixed an outside decomposition $\Pi$ of $D$, we can determine for
each cell $x$ in $D$, lying in one of the ribbons $\theta_k$, 
whether $x$ \emph{goes up} or \emph{goes right} in $\Pi$:  
\begin{enumerate}
\item[$\bullet$] It goes up if the cell immediately north of $x$ lies in the same ribbon $\theta_k$, or
if $x$ is the northeasternmost cell of $\theta_k$ and lies on the top perimeter of $D$.
\item[$\bullet$] It goes right if the cell immediately east of $x$ lies in the same ribbon $\theta_k$, 
or if $x$ is the northeasternmost cell of $\theta_k$ and lies on the right perimeter of $D$.
\end{enumerate}

A basic fact about outside decompositions $\Pi$ is that cells in the same diagonal
within $D$ will either all go up or all go right with respect to $\Pi$.  One can thus define
the \emph{cutting strip} $\theta(\Pi)$
for $\Pi$ to be the unique ribbon occupying the same nonempty diagonals
as $D$, such that the cell in a given diagonal goes up/right exactly as the cells of $D$ all do 
with respect to $\Pi$. Observe that each ribbon $\theta_k$ can be identified naturally with a subdiagram of the cutting strip $\theta(\Pi)$, and hence is uniquely determined
by the interval of contents $[p(\theta_k),q(\theta_k)]$ that its cells occupy.
In this way we can identify intervals $[p,q]$  with subribbons $\theta[p,q]$ of the
cutting strip $\theta(\Pi)$,
where we adopt the conventions that 
\begin{enumerate}
\item[$\bullet$] $\theta[q+1,q]$ represents the empty ribbon $\varnothing$,
having corresponding skew Schur function $s_\varnothing:=1$, and
\item[$\bullet$] $\theta[p,q]$ is \emph{undefined}
when $p > q+1$, and has corresponding skew Schur function $s_{\theta[p,q]}=0$.
\end{enumerate}

Using these conventions, define a new ribbon 
$$
\theta_i \# \theta_j:=\theta[p(\theta_j),q(\theta_i)]
$$ 
inside the cutting strip $\theta(\Pi)$.
Then the \emph{Hamel-Goulden determinant} formula asserts that

\begin{theorem}\cite{HamelGoulden}\label{Hamel-Goulden-formula} For any outside decomposition $\Pi=(\theta_1,\ldots,\theta_m)$ of a skew diagram $D$
$$
s_D = \det( s_{\theta_i \# \theta_j} )_{i,j=1}^{m}.
$$\end{theorem}

\begin{example} \rm \ \\
Consider the following skew diagram $D$, whose southwesternmost cell
is assumed to be $(1,1)$ with content $0$, and outside decomposition 
$\Pi=(\theta_1,\theta_2,\theta_3)$ where the cells in $\theta _i$ are labelled by $i$. Observe the associated  cutting strip $\theta(\Pi)$,
and the identification of the ribbons $\theta_k$ with intervals of contents within $\theta(\Pi)$:
$$
D=
\begin{array}{ccccc}
 & & &1&1\\
 &3&3&2&2\\
 &3&2&2& \\
3&3&2& & \\
3&2&2&    
\end{array},
\quad \quad
\theta(\Pi)=
\begin{array}{ccccc}
      &      &\times&\times&\times\\
      &\times&\times&      &      \\
      &\times&      &      &      \\
\times&\times&      &      &      \\
\times&      &      &      &
\end{array}
\quad\quad
\begin{aligned}
\theta_1&\leftrightarrow \theta[7,8] \\
\theta_2&\leftrightarrow \theta[1,7] \\
\theta_3&\leftrightarrow \theta[0,5].
\end{aligned}
$$
The associated Hamel-Goulden determinant is
$$
s_D
=\det
\left[
\begin{matrix}
s_{\theta[7,8]} & s_{\theta[1,8]} & s_{\theta[0,8]} \\
s_{\theta[7,7]} & s_{\theta[1,7]} & s_{\theta[0,7]} \\
s_{\theta[7,5]} & s_{\theta[1,5]} & s_{\theta[0,5]} \\
\end{matrix}
\right]
$$
$$
=\det
\left[
\begin{matrix}
  & &  & &\\
 &  
 s_{\begin{matrix} \times &\times  \\
    \end{matrix}} & &
  s_{\begin{matrix}     & &\times&\times&\times \\
                   &\times&\times&      &       \\
                   &\times&      &      &       \\
         \times&\times&      &      &       \\  
     \end{matrix}} & &
    s_{\begin{matrix}
                      &      &\times&\times&\times\\
                   &\times&\times&      &      \\
                 &\times&      &      &      \\
         \times&\times&      &      &      \\
          \times&      &      &      &      \\
       \end{matrix}} \\ 
 & & & &\\
 &
 s_{\begin{matrix} \times \end{matrix}} & &
  s_{\begin{matrix}     & &\times&\times \\
                       &\times&\times&       \\
                       &\times&      &       \\
              \times&\times&      &       \\
     \end{matrix}} & &
    s_{\begin{matrix}
            &      &\times&\times\\
            &\times&\times&      \\
            &\times&      &      \\
      \times&\times&      &      \\
      \times&      &      &      \\
       \end{matrix}} \\
  & & & &\\
  &
  0 & &
  s_{\begin{matrix} &  \times&\times&       \\
                    &         \times&      &       \\
                    \times&\times&      &       
     \end{matrix}} & &
    s_{\begin{matrix}
            &\times&\times&      \\
            &\times&      &      \\
      \times&\times&      &      \\
      \times&      &      &      
       \end{matrix}} \\
\end{matrix}
\right].
$$
\end{example}

There are two particular canonical outside decompositions of a connected skew diagram that will play an important role later.  

\begin{definition}\rm \ \\
\label{SE-NW-decompositions-definition}Given a connected skew diagram $D$, the \emph{southeast decomposition} is the 
following decomposition into ribbons, which is unique up to reordering. The first ribbon $\theta$ starts 
at the cell on the lower left, traverses the \emph{southeast} border of $D$, and ends at the cell on  
the upper right. Now consider $D$ with $\theta$ removed, which may decompose into several connected
component skew diagrams, and iterate the above procedure on each of these shapes in any order.
The \emph{northwest decomposition} is similarly defined, starting with a ribbon $\theta$ that
traverses the \emph{northwest} border of $D$.
\end{definition}

\begin{example} \rm \ \\
\label{SE-decomposition-example}For the following skew diagram $D$,  there are four ribbons 
$\theta_1,\theta_2,\theta_3,\theta_4$ in its southeast decomposition,
indicated by the numbers $1,2,3,4$ respectively:
$$
D=
\begin{matrix}
 & & & & & & & &4&4&3&3&3&1&1&1\\
 & & & & & & & &3&3&3&1&1&1& & \\
 & & & & & & &3&3&1&1&1& & & & \\
 & & & & & & &1&1&1& & & & & & \\
 & &2&2&2&1&1&1& & & & & & & & \\
 &2&2&1&1&1& & & & & & & & & & \\
 &1&1&1& & & & & & & & & & & & \\
1&1& & & & & & & & & & & & & & 
\end{matrix}.
$$
Here the first and largest ribbon $\theta=\theta_1=(2,3,3,3,3,3,3,3)$.
\end{example}

Note that for any connected skew diagram $D$,
both the southeast decomposition and the northwest decomposition
are outside decompositions of $D$, and hence give rise
to Hamel-Goulden determinants for $s_D$.  In both cases, the
associated cutting strip for this outside decomposition coincides with its first and largest ribbon $\theta$.

\section{The Littlewood-Richardson rule}\label{Lrrule}

  The Littlewood-Richardson rule gives the unique expansion of the skew Schur function
$s_{\lambda/\mu}$ into Schur functions $s_\nu$ for partitions $\nu$, and has many equivalent
versions.  We will use here a version suited to our purposes, which is known to be equivalent
to Zelevinsky's {\it picture} formulation \cite{Zelevinsky} of the rule.

\begin{definition} \rm \ \\
Given a skew diagram $D$, let its {\it row filling} $T_{\row}(D)$ be the function from
cells of $D$ to the integers which assigns to a cell its row index.

Say that a column-strict tableau $T$ is a {\it picture} for $D$ if 
\begin{enumerate}
\item[(i)]
the content of $T$ is the same as that of $T_{\row}(D)$, and 
\item[(ii)]  the map $f$ from cells of $D$ to cells of $T$, defined by sending
the $k^{th}$ cell from the {\it right} end of row $r$ of $D$ to the $k^{th}$ occurrence of 
the entry $r$ from the {\it left} in $T$, enjoys this additional property:
if a cell $x$ lies  lower in the same column of $D$ as some cell $x'$, then 
$f(x)$ lies in a lower row of $T$ than $f(x')$ (but not necessarily in the same column). 
\end{enumerate}
Denote by $\Pictures(D)$ the set of all column-strict tableaux that are
pictures for $D$.  Given a column-strict tableau $T$, let $\lambda(T)$ denote the partition
that gives its shape.

\end{definition}

\begin{theorem}(Littlewood-Richardson rule)
\begin{equation}
\label{LR-rule}
s_D = \sum_{T \in \Pictures(D)} s_{\lambda(T)}.
\end{equation}
\end{theorem}

\begin{example} \rm \ \\\label{LR-example}Consider the following skew diagram $D$, and its row filling $T_{\row}(D)$:
$$
D=\begin{matrix}
      &      &      &\times\\
      &      &      &\times\\
      &\times&\times&      \\
\times&\times&\times&      \\
\times&\times&      &      
\end{matrix},
\qquad
T_{\row}(D)=\begin{matrix}
 & & &1\\
 & & &2\\
 &3&3& \\
4&4&4& \\
5&5& &      
\end{matrix}.
\qquad
$$
Then one has 
$$
\begin{aligned}
&\Pictures(D)=\\
&\quad
\left\{
\begin{matrix}
1&3&4\\
2&4& \\
3&5& \\
4& & \\
5& & 
\end{matrix},\qquad
\begin{matrix}
1&3&4\\
2&4&5\\
3& & \\
4& & \\
5& & 
\end{matrix},\qquad
\begin{matrix}
1&3&4\\
2&4&5\\
3&5& \\
4& &
\end{matrix},\qquad
\begin{matrix}
1&3&3\\
2&4&4\\
4&5& \\
5& &
\end{matrix},\qquad
\begin{matrix}
1&3&3\\
2&4&4\\
4&5&5
\end{matrix}, \right. \\
&\qquad
\left.
\begin{matrix}
1&3&3&4\\
2&4& &\\
4&5& &\\
5& & &
\end{matrix},\quad 
\begin{matrix}
1&3&3&4\\
2&4&5&\\
4& & &\\
5& & &
\end{matrix},\quad 
\begin{matrix}
1&3&3&4\\
2&4&5&\\
4&5& &
\end{matrix},\quad 
\begin{matrix}
1&3&3&4\\
2&4&4&\\
5&5& &
\end{matrix},\quad 
\begin{matrix}
1&3&3&4\\
2&4&4&5\\
5& & &
\end{matrix} \right\}.
\end{aligned}
$$
Consequently the Littlewood-Richardson rule says
$$
\begin{aligned}
s_D &= s_{(3,2,2,1,1)} + s_{(3,3,1,1,1)} +2 s_{(3,3,2,1)} +s_{(3,3,3)}\\
     &\qquad +s_{(4,2,2,1)} + s_{(4,3,1,1)} + 2 s_{(4,3,2)} + s_{(4,4,1)}.
\end{aligned}
$$
\end{example}

\section{Reduction to connected diagrams\label{reduction-to-connected}}
We are now ready to state our key definition.

\begin{definition} \rm \ \\ 
Given two skew diagrams $D_1$ and $D_2$,
say that they are \emph{skew-equivalent}, denoted $D_1\sim D_2$, if $s_{D_1}=s_{D_2}$.
\end{definition}

The goal of this section is to understand two reductions:
\begin{enumerate}
\item[$A$.] Understanding all binomial syzygies among the skew Schur functions
is equivalent to understanding the equivalence relation $\sim$ on all skew diagrams, and
\item[$B$.] the latter is equivalent to understanding $\sim$ among {\it connected} skew diagrams.
\end{enumerate}

Both of these reductions will follow from some simple observations about the matrix 
$$
JT(\lambda/\mu):= (h_{\lambda_i - \mu_j -i + j})_{i,j=1}^{\ell(\lambda)},
$$
which appears in the Jacobi-Trudi determinant  \eqref{Jacobi-Trudi-formula} for a skew diagram $\lambda /\mu$.

\begin{proposition} 
\label{JT-diagonals-proposition}
Let $\lambda/\mu$ be a skew diagram with $\ell:=\ell(\lambda)$.
\begin{enumerate}[(i)]
\item
The largest subscript $k$ occurring on any nonzero entry $h_k$ in
the Jacobi-Trudi matrix $JT(\lambda/\mu)$ is 
$$
L :=\lambda_1+\ell-1  
$$
and this subscript occurs exactly once, on the $(1,\ell)$-entry $h_L$.
\item
The subscripts on the diagonal entries in $JT(\lambda/\mu)$ are
exactly the {row lengths}
$$
(r_1,\ldots,r_{\ell}):=(\lambda_1-\mu_1,\ldots,\lambda_{\ell}-\mu_{\ell})
$$
and the monomial $h_{r_1} \cdots h_{r_\ell}$ occurs in the determinant $s_D$
  \begin{enumerate}
     \item with coefficient $+1$, and
     \item as the monomial whose subscripts rearranged into weakly decreasing order
         give the smallest partition of $|\lambda/\mu|$ in dominance order among all nonzero monomials.
  \end{enumerate}
\item
The subscripts on the nonzero subdiagonal entries in $JT(\lambda/\mu)$ are
exactly one less than the {adjacent row overlap lengths}:
$$
(\lambda_2-\mu_1,\lambda_3-\mu_2,\ldots,\lambda_{\ell}-\mu_{\ell-1}).
$$

\end{enumerate}
\end{proposition}
\begin{proof}
Assertion (i) follows since  the subscripts appearing on nonzero entries in
$JT(\lambda/\mu)$ are of the form $\lambda_i-\mu_j-i+j$ with 
$$
\begin{aligned}
\lambda_i \leq \lambda_1, \
\mu_j     \geq 0, \
 i         \geq 1, \ 
j         \leq \ell
\end{aligned}
$$
so that  
$$
\lambda_i-\mu_j-i+j \,\, \leq \,\, \lambda_1-0-1+\ell \,\, =\,\, L.
$$
Furthermore, equality can occur only if $i=1$ and $j = \ell$.

For assertion (ii), expand the determinant of $JT(\lambda/\mu)$ as a signed
sum over of permutations in $\Sym_\ell$.  We claim that only the identity permutation
gives rise to the monomial $h_{r_1} \cdots h_{r_\ell}$.  This is because any other permutation
$\sigma$ can be obtained from the identity by a sequence of transpositions each increasing
the number of inversions, and it is straightforward to check that any such transpositions alters the
corresponding monomial so as to make its subscript sequence 
go \emph{strictly} upwards in the dominance order on partitions of $|\lambda/\mu|$.

Assertion (iii) is straightforward from the definitions, noting that
$\lambda_{i+1} - \mu_{i}$ is indeed the number of columns of overlap between
row $i$ and row $i+1$ in the skew diagram.
\end{proof}

\begin{corollary}
\label{reduction-to-connected-corollary} 
For a disconnected skew diagram $D= D_1 \oplus D_2$, one has the factorization $s_D=s_{D_1} s_{D_2}$. 
For a connected skew diagram $D$, the polynomial $s_D$ is irreducible in $\ZZ[h_1,h_2,\ldots]$.
\end{corollary}
\begin{proof}
The first assertion of the proposition is well-known, and follows, for example,
immediately from the definition 
\eqref{tableau-Schur-definition} of $s_D$ using tableaux.

For the second assertion on irreducibility\footnote{The proof 
of this irreducibility was incorrect in previous
versions of this paper, and is incorrect in the journal version ({\it Adv. Math.} (216) (2007),
118--152).  The authors are grateful to Marc van Leeuwen for bringing this to their attention.
A correct proof, obtained jointly with Farzin Barekat and to appear as a corrigenda
in Adv. Math., is supplied here.},
we will induct on the number $\ell:=\ell(D)$ of nonempty rows in the connected skew diagram $D$.
The base case $\ell=1$ is trivial, as then $s_D = h_{|D|}$ where $|D|$ is
the number of cells of $D$.

Thus, in the inductive step one may assume $\ell \geq 2$, and
assume for the sake of contradiction that $s_D$ is reducible.
Express $D=\lambda/\mu$ with $|\lambda|$ minimal, so that, 
in particular, $\ell(D)=\ell(\lambda)=:\ell$ and 
$\mu_\ell=0$.  Let $L=\lambda_1+\ell-1$, so that by 
Proposition~\ref{JT-diagonals-proposition}(i) 
the $\ell \times \ell$ Jacobi-Trudi matrix $J$ for $s_D$ expresses
\begin{equation}
\label{J-T-consequence}
s_D = s \cdot h_L + r,
\end{equation}
in which both $r, s$ involve only the variables $h_1,h_2,\ldots,h_{L-1}$. 

We claim that neither $r$ nor $s$ is the zero polynomial. 
For $r$, note that Proposition~\ref{JT-diagonals-proposition}(ii) implies that
 $r$ must contain the monomial $h_{r_1} \cdots h_{r_\ell}$
with coefficient $+1$ where $r_1, \ldots , r_\ell$ are the lengths of the rows of $\lambda /\mu$.
For $s$, note that $s$ is $(-1)^{\ell-1}$ times the determinant of the $(\ell -1)\times (\ell -1)$ complementary minor to $h_L$ in $J$, and  the complementary minor is the Jacobi-Trudi matrix for $s_{\hat{\lambda}/\hat{\mu}}=s_{\hat{D}}$, where $\hat{\lambda}=(\lambda _2, \lambda_3,\ldots , \lambda _\ell)$, 
$\hat{\mu} = (\mu _1 +1, \mu _2+ 2, \ldots, \mu_\ell +1)$.
Observe that $\hat{D}$ is obtained from $D$ by removing the northwesternmost
ribbon  from the northwest border of the connected skew diagram $D$ (in english notation).

Thus,  \eqref{J-T-consequence} shows that $s_D$ is {\it linear} as a polynomial in $h_L$.  Since we are assuming $s_D$ is reducible, 
this means $s_D$ must have at least one nontrivial irreducible factor, call it $f$, which is of degree zero in $h_L$.
This factor $f$ must therefore also divide $r$, and hence also divide $s$.

Denote by $J^A_B$ the submatrix obtained from $J$ by removing
its rows indexed by the subset $A$ and columns indexed by the subset $B$.
Then the {\it Lewis Carroll} or {\it Dodgson condensation} or 
{\it Desnanot-Jacobi adjoint matrix} identity \cite[Theorem 3.12]{Bressoud}
asserts that
\begin{equation}
\label{Desnanot-Jacobi}
\det J^{1,\ell}_{1,\ell} \cdot \det J = 
  \det J^{1}_{1} \cdot \det J^{\ell}_{\ell}
   -\det J^{1}_{\ell} \cdot \det J^{\ell}_{1}.
\end{equation}
Note that the left side of \eqref{Desnanot-Jacobi}
is divisible by $f$ since $\det J = s_D$, and the second term on the
right side of \eqref{Desnanot-Jacobi} is also divisible by $f$,
since $\det J^{1}_{\ell}$ is the same as the minor determinant appearing in $s$ in  \eqref{J-T-consequence}.
Therefore, the first term on the right of
\eqref{Desnanot-Jacobi} is divisible by $f$, implying
that one of its factors $\det J^{1}_{1}$ or $\det J^{\ell}_{\ell}$ 
must be divisible by $f$.  However, one can 
check that these last two determinants are 
the Jacobi-Trudi determinants for the skew diagrams $E, F$ obtained 
from $D$ by removing its first, last row, respectively.
Since $E, F$ are connected skew diagrams with fewer rows than $D$, both $s_E, s_F$ are
irreducible by the inductive hypothesis.  Hence either $f=s_E$ or
$f=s_F$.  But since $f$ divides $s$, its degree satisfies
$$
\deg(f) \leq \deg(s)   =|D|-(\lambda_1 + \ell -1)
$$
and this last quantity is strictly less than both
$$
\begin{aligned}
\deg(s_E) &= |D|-(\lambda_1-\mu_1), \text{ and }\\
\deg(s_F) &= |D|-(\lambda_\ell-\mu_\ell)
\end{aligned}
$$
since $\ell \geq 2$.  This contradicts having either $f=s_E$ or
$f=s_F$, ending the proof.
\end{proof}

We can now infer Reductions A and B from the beginning of the section.  Given a binomial syzygy
$$
c \,\, s_{D_1} s_{D_2} \cdots s_{D_m} - c' \,\, s_{D'_1} s_{D'_2} \cdots s_{D'_m} = 0
$$
among the skew Schur functions, with coefficients $c, c'$ in any ring, 
the first assertion of Corollary~\ref{reduction-to-connected-corollary} allows one to
rewrite this as $c \,\, s_{D} = c'\,\, s_{D'},$ where
$$
\begin{aligned}
D &:=D_1 \oplus D_2 \oplus \cdots \oplus D_m\\
D'&:=D'_1 \oplus D'_2 \oplus \cdots \oplus D'_m.
\end{aligned}
$$

\noindent
Proposition~\ref{JT-diagonals-proposition}(ii) implies the
unitriangular expansion
$$
s_D= h_{\rho} + \sum_{\mu: \,\, \mu >_{dom} \rho} c_\mu h_\mu
$$
in which  $\rho$ is the weakly decreasing rearrangement of the row lengths in $D$.
As $s_{D'}$ has a similar expansion, this 
forces $c=c'$ above, and hence $s_D = s_{D'}$.
That is, $D \sim D'$, achieving Reduction A.

For Reduction B, use the fact that $\Lambda = \ZZ[h_1,h_2,\ldots]$ is a unique
factorization domain, along with Corollary~\ref{reduction-to-connected-corollary}.

%%%%%%% Sufficiencies start here %%%%%%%%%%%

\section{Sufficient conditions\label{sufficiencies}}

The most basic skew-equivalence is the following well-known fact.

\begin{proposition}\cite[Exercise 7.56(a)]{Stanley1}\label{antipodal-equivalence}
If $D$ is a skew diagram then  $D \sim D^*$, where $D^*$ is the antipodal rotation of $D$.
\end{proposition}

Recently it was also proved that

\begin{theorem} \cite[Theorem 4.1]{BilleraThomasvanWilligenburg}\label{HDL_original}
Two ribbons $\beta$ and $\gamma$ satisfy $\beta\sim\gamma$ if and only if for some $k$
$$\beta _1\circ \cdots \circ\beta _k \sim \gamma _1 \circ \cdots \circ\gamma _k,$$where for each  $i$ either $\gamma _i = \beta _i$ or $\gamma _i = \beta _i ^*$. 
\end{theorem}

It transpires that there are several other
constructions and operations on skew diagrams that give rise to more skew-equivalences.

\subsection{Composition with ribbons}
\label{composition-section}

We now show that the notation for the diagrammatic operations 
$\alpha \circ D$ and $D \circ \beta$ defined in Section~\ref{ribbon-operations}
are consistent with algebraic operations on 
skew Schur functions $s_D$.  These operations then lead to nontrivial skew-equivalences.

We begin by reviewing the presentation of the ring $\Lambda$ of symmetric functions
by the generating set of ribbon Schur functions $s_\alpha$.  Let
$\QQ[z_\alpha]$ denote a polynomial algebra in infinitely many variables $z_\alpha$
indexed by all compositions $\alpha$.

\begin{proposition} \cite[Proposition 2.2]{BilleraThomasvanWilligenburg}
The algebra homomorphism
$$
\begin{matrix}
\QQ[z_\alpha] &\rightarrow &\Lambda \\
    z_\alpha  & \mapsto    &s_\alpha
\end{matrix}
$$
is a surjection, whose kernel is the ideal generated by the relations
\begin{equation}
\label{ribbon-relation}
z_\alpha z_\beta - ( z_{\alpha \cdot \beta} + z_{\alpha \odot \beta} ).\qed
\end{equation}
\end{proposition}

\begin{corollary}
\label{circ-well-defined}
For a fixed skew diagram $D$ the map 
$$
\begin{matrix}
\QQ[z_\alpha] &\overset{(-) \circ s_D}{\longrightarrow} &\Lambda \\
    z_\alpha  & \longmapsto    &s_{\alpha \circ D}
\end{matrix}
$$
descends to a well-defined algebra map $\Lambda \longrightarrow \Lambda$.
In other words, for any symmetric function $f$, one can
arbitrarily write $f$ as a polynomial in ribbon Schur functions 
$f=p(s_\alpha)$ and then set
$f \circ s_D:= p(s_{\alpha \circ D})$.
\end{corollary}

\begin{proof}
The relation \eqref{ribbon-relation} maps under  $(-) \circ s_D$ to
$$
\begin{aligned}
s_{\alpha \circ D} s_{\beta \circ D} 
    & -  \left( s_{(\alpha \cdot \beta) \circ D} 
         + s_{(\alpha \odot \beta) \circ D} \right) \\
= s_{\alpha \circ D} s_{\beta \circ D} 
    & -  \left( s_{(\alpha \circ D) \cdot (\beta \circ D)} 
         + s_{(\alpha \circ D) \odot (\beta \circ D)} \right)
\end{aligned}
$$
using Proposition~\ref{circ-distributivities}.  This last expression is zero by
Proposition~\ref{most-basic-syzygy}.
\end{proof}

We should point out that the notation $f \mapsto f \circ s_D$ has
already been used in \cite{BilleraThomasvanWilligenburg}  to denote the  \emph{plethysm}
or \emph{plethystic composition}, following one of the
standard references \cite{Macdonald}.  We will instead use the notation 
$f \mapsto f[s_\alpha]$ for plethysm, freeing the symbol $\circ$ for
use in the map $f \mapsto f \circ s_D$ defined in Corollary~\ref{circ-well-defined}.
Note that we are abusing notation by using $\circ$ both for the map $(-) \circ s_D$ on symmetric
functions, as well as the two diagrammatic operations $\alpha \circ D$
and $D \circ \beta$.  The previous corollary says that it is well-defined to set
\begin{equation}
\label{consistent-with-first}
s_\alpha \circ s_D = s_{\alpha \circ D}
\end{equation}
so that we are at least consistent with one of the diagrammatic operations.
The next result says that we are also consistent with the other.

\begin{proposition}
\label{consistent-with-second}
For any skew diagram $D$ and ribbon $\beta$
$$s_{D \circ \beta} = s_D \circ s_\beta.$$

\end{proposition}
\begin{proof}
Pick an outside decomposition $\Pi=(\theta_1,\ldots,\theta_m)$ of $D$, with cutting strip $\theta(\Pi)$, 
so that Theorem~\ref{Hamel-Goulden-formula} asserts $s_D = \det( s_{\theta_i \# \theta_j} )_{i,j=1}^{m}$.
It follows from the definition of $D\circ \beta$ and the definition of outside decomposition that
\begin{enumerate}
\item[$\bullet$]
$\Pi \circ \beta := (\theta_1 \circ \beta ,\ldots,\theta_m \circ \beta)$
gives an outside decomposition for $D \circ \beta$ and consequently from the definition of cutting strip
\item[$\bullet$]  that the cutting strip  satisfies the formula
$$
\theta( \Pi \circ \beta ) = \theta(\Pi) \circ \beta.
$$
\item[$\bullet$] Moreover, the relevant subribbons of this cutting strip satisfy the commutation
$$
( \theta_i \circ  \beta ) \# ( \theta_j \circ \beta ) = ( \theta_i \# \theta_j ) \circ \beta .
$$
\end{enumerate}
Consequently,
$$
\begin{aligned}
s_{D \circ \beta} 
  &=  \det\left[ s_{(\theta_i \circ \beta ) \# (\theta_j \circ \beta)} \right]^m_{i,j=1} \\
  &=  \det\left[ s_{(\theta_i \# \theta_j) \circ \beta } \right]^m_{i,j=1} \\
  &=  \det\left[ s_{\theta_i \# \theta_j} \right]_{i,j=1}^{m}  \circ s_\beta   \\
  &=   s_D \circ s_\beta,
\end{aligned}
$$where the third equality follows from Corollary~\ref{circ-well-defined}.
\end{proof}

We are now ready to state the first of two main ways to create new skew-equivalences from known ones.

\begin{theorem}\label{composition-equivalences}
Assume one has ribbons $\alpha, \alpha'$ and skew diagrams $D, D'$ satisfying
$\alpha \sim \alpha'$ and $D \sim D'$.  Then
\begin{enumerate}
\item[(i)]$\alpha \circ D \sim     \alpha' \circ D,$ 
\item[(ii)] $D \circ \alpha \sim D' \circ \alpha,$ 
\item[(iii)] $D \circ \alpha   \sim    D \circ \alpha',$ and
\item[(iv)]$\alpha \circ D \sim    \alpha \circ D^\ast.$
\end{enumerate}
\end{theorem}
\begin{proof}
Assertions (i) and (ii) both follow from the fact that if $E$ is any skew diagram, then
$D \sim D'$ means $s_D = s_{D'}$, and hence 
\begin{equation}
\label{mediating-equality}
s_D \circ s_E = s_{D'} \circ s_E.
\end{equation}
Note that
if $D, D'$ happen to be ribbons $\alpha, \alpha'$, then this gives the middle equality in
$$
s_{\alpha \circ E} = s_{\alpha} \circ s_E = s_{\alpha'} \circ s_E = s_{\alpha' \circ E},
$$
while \eqref{consistent-with-first} gives the outside equalities. 
This proves  $\alpha \circ E \sim \alpha' \circ E$,
and hence assertion (i).  Similarly, if $E$ happens to be a ribbon $\alpha$, then
\eqref{mediating-equality} again gives the middle equality in
$$
s_{D\circ \alpha} = s_D \circ s_\alpha = s_{D'} \circ s_\alpha = s_{D' \circ \alpha},
$$
and Proposition~\ref{consistent-with-second} gives the outside equalities.
This proves $D \circ \alpha \sim D' \circ \alpha$, and hence assertion (ii).

For assertion (iii), we deduce it first in the special
case where the skew diagram $D$ is a ribbon $\beta$.  It follows then
from Theorem~\ref{HDL_original}.  This characterization
asserts that $\alpha \sim \alpha'$
for two ribbons $\alpha, \alpha'$ if and only if there are expressions
\begin{equation}
\label{ribbon-skew-equivalence-characterization}
\begin{aligned}
\alpha&=\gamma_1 \circ \gamma_2 \circ \cdots \circ \gamma_r \\
\alpha'&=\delta_1 \circ \delta_2 \circ \cdots \circ \delta_r
\end{aligned}
\end{equation}
in which for each $i$ one has that $\gamma_i, \delta_i$ are ribbons 
with either $\gamma_i=\delta_i$ or $\gamma_i=\delta_i^\ast$.  Composing the
expressions in \eqref{ribbon-skew-equivalence-characterization}
with $\beta$ leads to similar such expressions for $\beta \circ \alpha, \beta \circ \alpha'$,
and hence $\beta \circ \alpha   \sim    \beta \circ \alpha'$. 

With this in hand, assertion (iii) for an arbitrary skew diagram $D$ is deduced
as follows.  Arbitrarily express $s_D=p(s_\beta)$ as a polynomial in various ribbon Schur functions 
$s_\beta$.  One then has the following string of equalities:
$$
s_{D \circ \alpha} 
 \overset{1}{=} s_D \circ s_\alpha 
 \overset{2}{=}  p(s_\beta) \circ s_\alpha 
 \overset{3}{=}  p(s_{\beta \circ \alpha})
 \overset{4}{=}  p(s_{\beta \circ \alpha'})
 \overset{5}{=}  p(s_\beta) \circ s_{\alpha'} 
 \overset{6}{=}  s_D \circ s_{\alpha'}
 \overset{7}{=}  s_{D \circ \alpha'}.
$$
Here the equalities $\overset{1}{=}$ and $\overset{7}{=}$ 
use Proposition~\ref{consistent-with-second},
the equalities $\overset{2}{=}$ and $\overset{6}{=}$ use the expression $s_D=p(s_\beta)$, 
the equalities $\overset{3}{=}$ and $\overset{5}{=}$ use Corollary~\ref{circ-well-defined},
and the equality $\overset{4}{=}$ uses the special case of (iii)
proven in the previous paragraph. Hence $D \circ \alpha \sim D \circ \alpha'$.

Assertion (iv) follows from assertion (i) and Proposition~\ref{antipodal-equivalence}:
$$
\alpha \circ D  \sim (\alpha \circ D)^\ast  = \alpha^\ast \circ D^\ast \sim \alpha \circ D^\ast.
$$\end{proof}

\begin{remark} \rm \ \\
Observe that Theorem~\ref{composition-equivalences} generalizes  \cite[Theorem 4.4 parts 1 and 2]{BilleraThomasvanWilligenburg}.  
\end{remark}

\begin{example} \rm \ \\
In general it is not true that $D \sim D'$ implies $\alpha \circ D \sim \alpha \circ D'$. For example, let $D=(4,3,2,1)/(1,1)$ and $D'=(4,3,2,1)/(2)$ and $\alpha =(2)$. Then $D\sim D'$ by Corollary~\ref{equivalent-conjugate-theorem}. However, $(8,7,6,5,3,2,1)/(5,5,4,1,1)=\alpha \circ D \not\sim \alpha \circ D'=(8,7,6,5,3,2,1)/(6,4,4,2)$ by Corollary~\ref{overlaps-corollary} below.
\end{example}
 
\begin{remark} \rm \ \\
It was observed in \cite[Proposition 3.4]{BilleraThomasvanWilligenburg} that
even though the $\circ$-composition and plethystic composition operations 
$$
\begin{aligned}
s_\alpha &\mapsto s_{\alpha \circ \beta} (=s_{\alpha} \circ {s_\beta})\\
s_\alpha &\mapsto s_{\alpha}[s_\beta]
\end{aligned}
$$
are {\it not} the same, they {\it do} coincide when one {\it sums/averages} over
all compositions $\alpha$ of a fixed size $n$:
$$
\left(\sum_{\alpha \vDash n} s_\alpha \right) \circ s_\beta
  = \sum_{\alpha \vDash n} s_{\alpha \circ \beta}
  = (s_\beta)^n 
  = (s_1^n)[s_\beta] 
  = \left( \sum_{\alpha \vDash n} s_\alpha \right) [s_\beta]
$$
in which the second and fourth equalities comes from iterating Proposition \ref{most-basic-syzygy}.  The
same holds replacing $s_\beta$ by $s_D$ for any skew diagram $D$, with the same proof:
$$
\left(\sum_{\alpha \vDash n} s_\alpha \right) \circ s_D
  = \sum_{\alpha \vDash n} s_{\alpha \circ D}
  = (s_D)^n
  = (s_1^n)[s_D]
  = \left( \sum_{\alpha \vDash n} s_\alpha \right) [s_D].
$$
\end{remark}

\begin{remark} \rm \ \\
\label{smallest-nonribbon-circ-example}Even though we have not defined a skew diagram $D_1 \circ D_2$ when $D_1, D_2$ are {\it both} non-ribbon
skew diagrams, the symmetric function $s_{D_1} \circ s_{D_2}$ is still well-defined,
via Corollary~\ref{circ-well-defined}.  One might ask whether there exists a
skew diagram $D$ playing the role of $D_1 \circ D_2$, that is, with
$s_D = s_{D_1} \circ s_{D_2}$.
Curiously and suggestively, computer calculations show that
this seems to be the case in the smallest example, in which $D_1,D_2$
are both the $2 \times 2$ rectangular Ferrers diagram  $\lambda=(2,2)$:
$$
s_{(2,2)} \circ s_{(2,2)} = s_{(1)} s_{(5,5,4,4,2)/(3,1,1)}.
$$
In other words, $(2,2) \circ (2,2)$ cannot be chosen to be a connected skew diagram, but
rather should be defined as the direct sum of a single cell with
$$
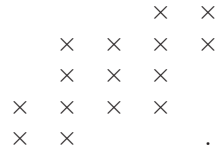
\begin{matrix}
      &      &      &\times&\times \\
      &\times&\times&\times&\times\\
      &\times&\times&\times&\\
\times&\times&\times&\times&\\
\times&\times&      &      &.\\
\end{matrix}
$$
This is somewhat remarkable, and  suggests a further avenue
of investigation for skew-equivalences; see Section \ref{Conclusion} below.
\end{remark}

%%%%%%% AMALGAMATION %%%%%%%%%%

\subsection{Amalgamation and amalgamated composition of ribbons}
\label{amalgamation-section}

In this section we introduce an operation $\alpha \circ_\omega D$ for certain skew diagrams $D$ and ribbons $\omega$, 
which we will call the {\it amalgamated composition} of $\alpha$ and
$D$ with respect to $\omega$.  It is analogous to the operation $\alpha \circ \beta$ on ribbons $\alpha, \beta$ and allows us to identify more 
skew diagrams that are skew-equivalent.

\begin{definition} \rm \ \\ 
Given a skew diagram $D$ and a nonempty ribbon $\omega$, 
say that $\omega$ {\it protrudes from the top (resp. bottom) of $D$} if  the restriction of $D$ to its $|\omega|$ northeasternmost 
(resp. southwesternmost) diagonals is the ribbon $\omega$ and the restriction of $D$ to its $|\omega|+1$ northeasternmost 
(resp. southwesternmost) diagonals is also a ribbon.

Given two skew diagrams $D_1, D_2$ and a nonempty ribbon $\omega$ protruding from  the top of $D_1$ and the
bottom of $D_2$, the \emph{amalgamation of $D_1$ and $D_2$ along $\omega$}, 
denoted $D_1 \amalg_\omega D_2$,
is the new skew diagram obtained from the disjoint union $D_1 \oplus D_2$ by identifying 
the copy of $\omega$ in the northeast of $D_1$ with 
the copy of $\omega$ in the southwest of $D_2$.

\end{definition}

\begin{example} \rm \ \\
Consider the skew diagram
$$
D=
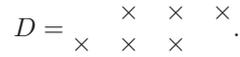
\begin{matrix}
       & \times & \times & \times \\
\times & \times & \times & 
\end{matrix}.
$$
Then $D$ has $\omega= \times$ protruding from the top and bottom.
 Furthermore,
$$
D \amalg_{\omega} D =
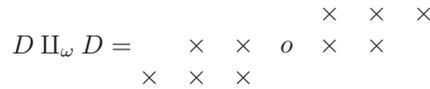
\begin{matrix}
       &       &       &        & \times & \times & \times \\
       & \times& \times& o & \times & \times &  \\
\times & \times& \times&        &        &        &  \\
\end{matrix}$$
and the copies of $\omega$  that have been amalgamated are
indicated with the letter $o$.
\end{example}

\begin{definition} \rm \ \\
When $\omega$ protrudes from the top of $D_1$ and bottom of $D_2$, one can
form the {\it outer (resp. inner) projection of $D_1$ onto $D_2$ with respect to $\omega$}.
This is a new diagram in the plane, not necessarily skew, 
obtained from the disjoint union $D_1 \oplus D_2$ by translating
$D_2$ and $D_1$ until the two copies of $\omega$ in $D_1, D_2$ are adjacent and
occupy the same set of diagonals, and the copy of $\omega$ in $D_1$ is immediately northwest (resp. southeast) of $\omega$ in $D_2$.

Observe that at most one of these two projections can be a skew diagram (and possibly neither one is).
When one of them is a skew diagram, denote it by $D_1 \cdot_{\omega} D_2$, and say that
{\it $D_1 \cdot_\omega D_2$ is defined} in this case. 
\end{definition}

\begin{example} \rm \ \\
Let  $D, \omega$ be as in the previous example.
Then the outer and inner projections of $D$ onto $D$ with respect to  $\omega$ are
$$
\begin{matrix}
       & \times& \times& o      &  &\times & \times & \times  \\
\times & \times& \times&   &   o   &\times & \times & 
\end{matrix} \qquad \qquad
\begin{matrix}
       &       &        & \times & \times & \times \\
       &       & o      & \times & \times & \\
       & \times& \times & o      &        & \\
\times & \times& \times &        &        & 
\end{matrix} = D \cdot_{\omega} D
$$
and only the latter is a skew diagram.
\end{example}

\begin{definition} \rm \ \\
Given a skew diagram $D$, and $\omega$ a ribbon protruding from both the top and bottom of $D$,
one can define 
$$
D^{\amalg_\omega n} 
  = \underbrace{D \amalg_\omega D \amalg_\omega \cdots \amalg_\omega D}_{n\text{ factors}} 
 := ((D \amalg_\omega D) \amalg_\omega D) \amalg_\omega \cdots \amalg_\omega D.
$$
If one assumes that $D \cdot_\omega D$ is also defined then by considering the northeasternmost copy of $D$ in $D^{\amalg_\omega m}$ and the southwesternmost copy of $D$ in $D^{\amalg_\omega n}$ 
for any positive integers $m,n$, we have
$(D^{\amalg_\omega m}) \cdot_\omega (D^{\amalg_\omega n})$ 
is also defined.
Under this assumption, for any ribbon $\alpha=(\alpha_1,\ldots,\alpha_\ell)$, define
the {\it amalgamated composition} of $\alpha$ and $D$ with respect to $\omega$ to be
the diagram
\begin{equation}
\alpha \circ_{\omega} D := 
( D^{\amalg_{\omega} \alpha_1} )
   \cdot_{\omega} \ldots \cdot_{\omega} 
     ( D^{\amalg_\omega \alpha_\ell}).
\label{composition_wrt_rho}\end{equation}
\end{definition}

\begin{example} \rm \ \\
Let $D, \omega$ be as in the previous example.  Then we saw earlier that
$D \cdot_\omega D$ is defined.  Consider the ribbon
$$
\alpha=(\alpha_1,\alpha_2,\alpha_3)=(2,1,3)=
\begin{matrix}
      &\times&\times&\times\\
      &\times&      &\\
\times&\times&      &.\\
\end{matrix}
$$
Then one has
$$
\begin{aligned}
\alpha \circ_\omega D 
&= (D^{\amalg_\omega 2}) \cdot_\omega 
        (D^{\amalg_\omega 1}) \cdot_\omega
          (D^{\amalg_\omega 3})\\
&=
\begin{matrix}
 & & & & & & & & & & & & & &3&3&3\\
 & & & & & & & & & & &3&3&3&3&3& \\
 & & & & & & & &3&3&3&3&3& & & & \\
 & & & & & & &3&3&3& & & & & & & \\
 & & & & & &2&2&2& & & & & & & & \\
 & & & & &2&2&2& & & & & & & & & \\
 & & & &1&1&1& & & & & & & & & & \\
 &1&1&1&1&1& & & & & & & & & & & \\
1&1&1& & & & & & & & & & & & & &
\end{matrix}.
\end{aligned}
$$
\end{example}

We now show that the operation $\alpha \circ_\omega D$ and the
$\alpha \circ D$ operation defined in Section~\ref{composition-section} 
associate with each other in a natural way.

\begin{proposition}
\label{amalgamation-associativities}
When $\alpha, \beta,\omega$ are ribbons and $D$ is a skew diagram such that
the appropriate operations are well-defined, one has
$$ (\alpha \circ \beta) \circ_\omega D = \alpha \circ_\omega (\beta \circ_\omega D).$$
\end{proposition}

\begin{proof}
This follows from the definitions since
\begin{eqnarray*}
(\alpha \circ \beta) \circ _\omega D
 &=&(D^{\amalg_{\omega} \beta _1}
   \cdot_{\omega} \cdots \cdot_{\omega} 
     D^{\amalg_\omega \beta _m})^{\amalg_{\omega} \alpha_1} \cdot_{\omega} \cdots \cdot_{\omega} \\
     &&(D^{\amalg_{\omega} \beta _1}
   \cdot_{\omega} \cdots \cdot_{\omega} 
     D^{\amalg_\omega \beta _m})^{\amalg_{\omega} \alpha_\ell}\\
     &=& \alpha \circ _\omega (\beta \circ _\omega D).
\end{eqnarray*}
\end{proof}

We now wish to interpret the diagrammatic operation $\alpha \circ_\omega D$ in
terms of an algebraic operation, for certain skew diagrams $D$ and ribbons $\omega$.

\begin{definition} \rm \ \\
Suppose that $D$ is a skew diagram and $\omega$ a ribbon protruding from the top and
bottom of $D$, so that $D^{\amalg_\omega r}$ is defined for all positive integers $r$.
Define a map of sets
$$
\begin{matrix}
\Lambda & \overset{(-) \circ_\omega s_D}{\longrightarrow}& \Lambda\\
    f   &      \longmapsto                               & f \circ_\omega s_D
\end{matrix}
$$
as the composite of two maps $\Lambda \rightarrow \Lambda[t] \rightarrow \Lambda$,
which we now describe.

Thinking of $\Lambda$ as the polynomial algebra $\ZZ[h_1,h_2,\ldots]$, we can
temporarily grade
$\Lambda$ and $\Lambda[t]$ by setting $\deg(t)=\deg(h_r)=1$ for all $r$. Note that this is \emph{not} the usual grading on
$\Lambda$, in which $\deg(h_r)=r$, and for which skew Schur functions $s_D$ are homogeneous.
In fact, $s_D$ will generally be {\it inhomogeneous} with respect to this temporary grading.
The first map $\Lambda \rightarrow \Lambda[t]$ simply homogenizes a polynomial in the $h_r$s
with respect to this grading, using the variable $t$ as the homogenization variable.  

The second map is defined by
$$
\begin{matrix}
\Lambda[t]& \longrightarrow &\Lambda \\
     h_r  & \longmapsto     &s_{D^{\amalg_\omega r}} \\
      t   & \longmapsto     &s_\omega
\end{matrix}
$$
Note that this composite map is not a ring homomorphism, nor even a map of
$\ZZ$-modules, because these properties fail for the homogenization map  $\Lambda \rightarrow \Lambda[t]$.
\end{definition}

Before we state the next theorem we need some hypotheses.

\begin{hypotheses}\label{wseparation}\rm \ \\
Suppose that $D$ is a connected skew diagram and $\omega$ is a ribbon protruding from the top and bottom of $D$. We assume that $D$ and $\omega$ satisfy the following conditions:
\begin{enumerate}
\item[(i)] $D \cdot_\omega D$ is defined,
\item[(ii)] the
two copies of $\omega$ protruding from the top and bottom of $D$ are separated by at least one
diagonal, that is, there is a nonempty diagonal in $D$ intersecting neither copy of $\omega$.
\end{enumerate}
\end{hypotheses}

\begin{theorem}
\label{amalgamated-composition-theorem}
Let $D$ be a connected skew diagram, and $\omega$ a ribbon satisfying Hypotheses~\ref{wseparation}. Then for any ribbon $\alpha$ one has
$$
s_{\alpha \circ_\omega D} = s_\alpha \circ_\omega s_D.
$$
\end{theorem}

\begin{remark} \rm \ \\
In Theorem~\ref{amalgamated-composition-theorem},
some hypothesis about separating the two copies of $\omega$ within $D$
is needed, as shown by the following example.  Let $\alpha$ be the ribbon $(1,1,1)$, let $D$ be the ribbon
$(1,1)$, and $\omega$ the single cell $(1)$.  In other words, let $\alpha, D, \omega$, respectively,
be diagrams that consist of a single column, of sizes $3,2,1$, respectively.

Then $\omega$ protrudes from the top and bottom of $D$, and one can check that
$$
D \cdot_\omega D= 
\begin{matrix} 
\times & \times \\
\times & \times
\end{matrix}
\qquad \text{ and }\qquad
\alpha \circ_\omega D=
\begin{matrix} 
\times & \times & \times\\
\times & \times & \times
\end{matrix}
$$
are defined.  However, the two copies of $\omega$ within $D$ occupy adjacent
diagonals, so that they fail the separation hypothesis in the theorem.
Correspondingly, one finds that
$$
\begin{aligned}
s_\alpha \circ_\omega s_D 
&= 
\det \left[ \begin{matrix}
h_1 & h_2 & h_3 \\
1   & h_1 & h_2 \\
0   & 1   & h_1
\end{matrix} \right] \circ_\omega s_D \\
&=
\det \left[ \begin{matrix}
s_D & s_{D \amalg_\omega D} & s_{D \amalg_\omega D \amalg_\omega D} \\
s_\omega   & s_D  & s_{D \amalg_\omega D} \\
0   & s_\omega   & s_D
\end{matrix} \right] \\
& = s_{\begin{matrix} \times & \times & \times \\ \,\,
                         \times & \times & \times \end{matrix}}
    - \,\, s_{\begin{matrix} \times & \times \\
                         \times & \times \\
                          \times & \times \end{matrix}} \\
& \neq s_{\alpha \circ_\omega D}.
\end{aligned}
$$
\end{remark}

\begin{proof} (of Theorem~\ref{amalgamated-composition-theorem})
We induct on the number of rows $k$ in the ribbon $\alpha$.  In the base case $k=1$,  by Equation~\eqref{Jacobi-Trudi-formula} one
has $s_\alpha= h_r$ for some $r$, and the assertion is trivial.

For the inductive step, let
$$
\begin{aligned}
\alpha &= (\alpha_1,\alpha_2,\alpha_3,\ldots,\alpha_k) \\
\bar{\alpha} &=(\alpha_2,\alpha_3,\ldots,\alpha_k) \\
\hat{\alpha} &=(\alpha_1+\alpha_2,\alpha_3,\ldots,\alpha_k).
\end{aligned}
$$
Then expanding the Jacobi-Trudi determinant for $s_\alpha$ along its last row gives
$$
s_\alpha = h_{\alpha_1} s_{\bar{\alpha}} - 1 \cdot s_{\hat{\alpha}}
$$
and hence that
\begin{equation}
\label{J-T-recursion}
\begin{aligned}
s_\alpha \circ_\omega s_D 
  &= (h_{\alpha_1} s_{\bar{\alpha}} - 1 \cdot s_{\hat{\alpha}}) \circ_\omega s_D \\
  &= s_{D^{\amalg_\omega \alpha_1}} (s_{\bar{\alpha}} \circ_\omega s_D) 
        - s_\omega (s_{\hat{\alpha}} \circ_\omega s_D) \\
  &= s_{D^{\amalg_\omega \alpha_1}} \, s_{\bar{\alpha} \circ_\omega D} 
        - s_\omega \, s_{\hat{\alpha} \circ_\omega D}
\end{aligned}
\end{equation}
where the last equality uses the inductive hypothesis.

We wish to compare this last expression with an expansion for a certain Hamel-Goulden
determinant computing $s_{\alpha \circ_\omega D}$.  Note that the two copies
of $\omega$ lying in the top and bottom of $D$
are subribbons of the longest ribbon in the southeast decomposition $\Pi$ of $D$, 
namely the cutting strip $\theta:=\theta (\Pi)$.
More generally, the two copies of $\omega$ in any diagram $D^{\amalg_\omega r}$ 
are subribbons of the longest ribbon in its southeast decomposition, namely
$\theta^{\amalg_\omega r}$.
One can then collate these southeast decompositions for $D^{\amalg_\omega \alpha_i}$ to produce
an outside decomposition $(\theta_1,\ldots,\theta_n)$ for
$$
\alpha \circ_\omega D \
 = D^{\amalg_\omega \alpha_1} \cdot_\omega \ldots \cdot_\omega D^{\amalg_\omega \alpha_k}
$$
in which the ribbons come in $k$ different blocks, with those in the 
$j^{th}$ block comprising the subdiagram $D^{\amalg_\omega \alpha_j}$.
Furthermore, because of the separation hypothesis about the two copies of $\omega$ in $D$,
ribbons in different blocks will almost never share any nonempty diagonals, as this will
only happen for the longest ribbon in two adjacent blocks.  For notational purposes below,
let $m$ be the number of ribbons in the first block, and index the longest ribbons
in the first and second blocks as $\theta_m$ and $\theta_{m+1}$.

Let $A$ be the Hamel-Goulden matrix for this outside decomposition of $\alpha \circ_\omega D$.
We will do a generalized Laplace expansion \cite[\S 1.8]{Shilov} of its determinant along the first $m$ rows.
Given subsets $R, C$ of $[n]:=\{1,2,\ldots,n\}$, let $A_{R,C}$ be the submatrix of $A$ having
rows and columns indexed by $R$ and $C$ respectively.  Then the generalized Laplace expansion says that
$$
\det A 
=\sum_{\substack{C \subset [n] \\ |C|=m}} \epsilon_C \,\, \det(A_{[m],C})  \det(A_{[m+1,n],[n] \backslash C})
$$
where $\epsilon_C = \pm 1$ is the sign of the permutation which sorts the concatenation of
$C$ and $[n] \backslash C$, both written in increasing order, to the sequence $1,2,\ldots,n$.

The foregoing observations about separation of diagonals imply that $A_{[m+1,n],[n] \backslash C}$ will
have a zero column (and hence vanishing determinant) unless the $m$-element subset $C$ is chosen
to contain all the columns $1,2,\ldots,m-1$, so that for some $j \in [m,n]$, one has
$C=[m-1] \cup \{j\}$ and hence $\epsilon_C = (-1)^{j-m}$.  Thus
$$
\begin{aligned}
s_{\alpha \circ_\omega D} 
  &=\sum_{j=m}^n (-1)^{j-m} \det(A_{[m],[m-1] \cup \{j\}}) \det(A_{[m+1,n],[m,n] \backslash \{j\}})\\
  &=s_{D^{\amalg_\omega \alpha_1}} \, s_{\bar{\alpha} \circ_\omega D} \\
  & \qquad   + \sum_{j=m+1}^n (-1)^{j-m} \det(A_{[m],[m-1] \cup \{j\}}) \cdot s_\omega \cdot
                                    \det(A_{[m+2,n],[m+1,n] \backslash \{j\}})
\end{aligned}
$$
where the last equality uses the fact that the first column of $A_{[m+1,n],[m,n]}$ contains only one
nonzero entry, namely $A_{m+1,m}= s_\omega$.  Comparing this with Equation \eqref{J-T-recursion},
it only remains to show that 
\begin{equation}
\label{identify-summation}
s_{\hat{\alpha} \circ_\omega D} = -
\sum_{j=m+1}^n (-1)^{j-m} \det(A_{[m],[m-1] \cup \{j\}}) \det(A_{[m+2,n],[m+1,n] \backslash \{j\}}).
\end{equation}
To see this, note that we can obtain an outside decomposition
of $\hat{\alpha} \circ_\omega D$ by starting with the outside decomposition
$(\theta_1,\ldots,\theta_n)$ for $\alpha \circ_\omega D$ used above, 
and replacing the two ribbons $\theta_m, \theta_{m+1}$ with a single ribbon 
$\theta_m \amalg_\omega \theta_{m+1}  = \theta^{\amalg_\omega \alpha_1+\alpha_2} .$
Now expand the corresponding $(n-1) \times (n-1)$ Hamel-Goulden determinant for 
$s_{\hat{\alpha} \circ_\omega D}$ along its first $m$ rows, and one obtains 
\eqref{identify-summation}.
\end{proof}

We are now ready to state our second key way to create new skew-equivalences from known ones.

\begin{theorem}
\label{RibbComp}
Let $\alpha, \alpha'$ be ribbons with $\alpha \sim \alpha'$, and
assume that $D, \omega$ satisfy Hypotheses~\ref{wseparation}.
Then one has the following skew-equivalences:
$$
\alpha' \circ_{\omega} D
\,\, \sim \,\,
\alpha \circ_{\omega} D 
\,\, \sim \,\,
\alpha \circ_{\omega^\ast} D^\ast
.$$ 
\end{theorem}
\begin{proof}
Both skew-equivalences are immediate from Theorem~\ref{amalgamated-composition-theorem}.
For the second, note that $(D^\ast)^{\amalg_{\omega^\ast} r} = (D^{\amalg_{\omega} r})^\ast$
for all $r$, so that the maps 
$$
\begin{aligned}
\Lambda &\overset{(-)\circ_\omega s_D}{\rightarrow} \Lambda \\
\Lambda &\overset{(-)\circ_{\omega^\ast} s_{D^\ast}}{\rightarrow} \Lambda
\end{aligned}
$$are the same.
\end{proof}

\begin{remark} \rm \ \\
Theorem~\ref{RibbComp} is analogous to  \cite[Theorem 4.4 parts 1 and 2]{BilleraThomasvanWilligenburg}.  
\end{remark}

\begin{theorem}\label{HDL_omega}
Let $\{ \beta _i \}_{i=1}^k, \{ \gamma _i \}_{i=1}^k$ be ribbons, 
and for each $i$ either $\gamma _i = \beta _i$ or $\gamma _i = \beta _i ^*$.
If the skew diagrams $D$, $\omega$ satisfy Hypotheses~\ref{wseparation}, then 
$$
\begin{aligned}
&\gamma _1 \circ _\omega \gamma _2 \circ _\omega \ldots \circ _\omega\gamma _k \circ _\omega D \\
&\quad \sim
\beta _1 \circ _\omega \beta _2 \circ _\omega \ldots \circ _\omega\beta _k \circ _\omega D \\
&\qquad \sim
\beta _1 \circ _{\omega ^*}\beta _2 \circ _{\omega ^*}\ldots \circ _{\omega ^*}\beta _k \circ _{\omega ^*} D^* 
\end{aligned}
$$
where all the operations $\circ _{\omega}$ or $\circ _{\omega^*}$ are performed from right to left.
\end{theorem}

\begin{proof}
By Theorem~\ref{HDL_original} and Theorem~\ref{RibbComp}
 we know
$$
(\gamma _1 \circ \ldots\circ\gamma _k) \circ _\omega D
\sim (\beta _1\circ \ldots \circ\beta _k)\circ _\omega D  
\sim  (\beta _1\circ \ldots \circ\beta _k)\circ _{\omega ^*} D^*.
$$
From \cite[Proposition 3.3]{BilleraThomasvanWilligenburg} we know $\circ$ is associative, 
and by applying Proposition~\ref{amalgamation-associativities} 
repeatedly $k-1$ times the result follows.
\end{proof}

\begin{remark}\rm \ \\
Theorem~\ref{HDL_omega} is analogous to the reverse 
direction of \cite[Theorem 4.1]{BilleraThomasvanWilligenburg}. 
\end{remark}

%%%%%%%%% RIBBON STAIRCASES %%%%%%%%%%

\subsection{Conjugation and ribbon staircases\label{ribbon-staircases}}
Recall from Definition~\ref{SE-NW-decompositions-definition} the \emph{southeast decomposition} and 
\emph{northwest decomposition} of a connected skew diagram.
When either of these decompositions takes on a very special form, we will show that it gives
rise to a nontrivial skew-equivalence, and in some cases to a skew-equivalence of the
form $D \sim D^t$.

\begin{definition}\rm\ \\
Let $\alpha =(\alpha _1, \ldots, \alpha _k)$ and $\beta = (\beta _1,\ldots,\beta _\ell)$ 
be ribbons.  For an integer $m \geq 1$, say that the {\it $m$-intersection} $\alpha \cap_m \beta$ {\it exists} 
if there is a ribbon $\omega=(\omega_1,\ldots,\omega_m)$
with $m$ rows protruding from the top of $\alpha$ and the bottom of $\beta$ for which $\omega_1=\beta_1$ and $\omega_m=\alpha_k$;
when $m=1$, we set $\omega_1:=\min\{\alpha_k, \beta_1\}$.  In this case,
define the {\it $m$-intersection} $\alpha \cap_m \beta$  and the {\it $m$-union} $\alpha \cup_m \beta$ to be
$$
\begin{aligned}
\alpha \cap_m \beta&:=\omega \\
\alpha \cup_m \beta&:=\alpha \amalg_\omega \beta .
\end{aligned}
$$
\noindent
If $\alpha \cup _m \beta = \alpha$ or $ \beta$ (resp. or $\alpha \cap _m \beta =\alpha$ or $\beta$) 
then we say the $m$-union (resp. $m$-intersection) is \emph{trivial}.
If $\alpha$ is a ribbon such that
$\alpha \cap _m\alpha$ exists and is nontrivial then
\begin{equation*} 
\varepsilon ^k _m(\alpha) :=
\underbrace{\alpha \cup _m \alpha \cup _m \ldots \cup _m \alpha}_{k\text{ factors}} 
\end{equation*}
is the \emph{ribbon staircase} of height $k$ and depth $m$ generated by $\alpha$.

%If $\alpha=\varepsilon ^\ell_n(\beta)$ with $\ell>1$ and $n<\ell(\alpha)$ for some ribbon $\beta$, 
%then we say $\alpha$ is \emph{ribbon decomposable}; otherwise $\alpha$ is \emph{ribbon indecomposable}.

\end{definition}

\begin{example}\rm \ \\
Let $\alpha$ be the ribbon $(2,3)$.
Then 
$$
\varepsilon ^3 _1(\alpha) = \varepsilon ^3 _1 \left(
\begin{matrix}
      &\times&\times&\times\\
\times&\times&      &
\end{matrix} \right)
=
\begin{matrix}
 & & & & &\times&\times&\times\\
 & & &\times&\times&\times\\
 &\times&\times&\times\\
\times&\times
\end{matrix}.
$$
\end{example}

\begin{definition} \rm\ \\
Say that a skew diagram $D$ has a {\it southeast ribbon staircase decomposition} 
if there exists an $m < \ell(\alpha)$ and a ribbon $\alpha$
such that all ribbons in the southeast decomposition of $D$ are of the form
$\alpha \cap_m \alpha$ or $\varepsilon^p_m (\alpha)$ for various integers $p\geq 1$.  

In this situation, let $k$ be the maximum value of $p$ occurring among the $\varepsilon^p_m (\alpha)$ above, so that
the largest ribbon $\theta$ equals $\varepsilon ^k _m(\alpha)$.  We will
think of $\theta$ as containing $k$ copies of $\alpha$, numbered
$1,2,\ldots,k$ from southwest to northeast. We now wish to define the {\it nesting} $\mathcal{N}$ associated to this decomposition.
The nesting $\mathcal{N}$ is a word of length $k-1$ using as letters the four symbols,
{\it dot} ``$.$'', {\it left parenthesis} ``$($'', {\it right parenthesis} ``$)$'' and
{\it vertical slash} ``$|$''.  Considering the ribbons in the southeast decomposition of $D$,
\begin{enumerate}
\item[$\bullet$]
a ribbon of the form $\varepsilon^p_m(\alpha)$ creates a pair of left and right parentheses
in positions $i$ and $j$ if the ribbon occupies the same diagonals
as the copies of $\alpha$ in $\theta$
numbered $i+1,i+2,\ldots,j-1,j$,   while
\item[$\bullet$]
a ribbon of the form $\alpha \cap_m \alpha$ creates a vertical slash in position $i$
if it occupies the same diagonals as the intersection of the $i,i+1$ copies
of $\alpha$ in $\theta$, and
\item[$\bullet$] all other letters in $\mathcal{N}$ are dots.
\end{enumerate}
With this notation, say that $D=(\varepsilon ^k _m(\alpha), \mathcal{N})_{se}$.
Analogously define the notation
$D=(\varepsilon ^k _m(\alpha), \mathcal{N})_{nw}$
using the northwest decomposition.

Lastly, given a nesting $\mathcal{N}$, denote the \emph{reverse nesting},  which is the reverse of the
word $\mathcal{N}$, by $\mathcal{N}^\ast$.

Observe that a nesting is well-defined as if we wanted to place two different parentheses, or a parenthesis or a slash, in a given position then this would imply either that we did not have a ribbon staircase decomposition or that we did not have a skew diagram.
\end{definition}

\begin{example}
\label{ribbon-staircase-example} \rm \ \\
 Recall the southeast decomposition of a skew diagram $D$ from Example~\ref{SE-decomposition-example}:
$$
D=
\begin{matrix}
 & & & & & & & &4&4&3&3&3&1&1&1\\
 & & & & & & & &3&3&3&1&1&1& & \\
 & & & & & & &3&3&1&1&1& & & & \\
 & & & & & & &1&1&1& & & & & & \\
 & &2&2&2&1&1&1& & & & & & & & \\
 &2&2&1&1&1& & & & & & & & & & \\
 &1&1&1& & & & & & & & & & & & \\
1&1& & & & & & & & & & & & & &
\end{matrix}.
$$
This is a southeast ribbon staircase decomposition, in which 
$\alpha=(2,3), m=1, k=7$ and 
$$
\begin{matrix}
\mathcal{N} &=& ( & ) & . & ( & | & ) \\
            & & 1 & 2 & 3 & 4 & 5 & 6
\end{matrix}
$$
that is, $D=(\varepsilon^7_1(\alpha), \mathcal{N} )_{se}$.
Here 
$$
\begin{matrix}
\mathcal{N}^\ast &=& ( & | & ) & . & ( & ) \\
                 & & 1 & 2 & 3 & 4 & 5 & 6
\end{matrix}
$$
and 
$D'=(\varepsilon^7_1(\alpha), \mathcal{N}^\ast )_{se}$ is the following skew diagram:
$$
D'=
\begin{matrix}
 & & & & & & & & & &2&2&2&1&1&1\\
 & & & & & & & & &2&2&1&1&1& & \\
 & & & & & & & & &1&1&1& & & & \\
 & &4&4&3&3&3&1&1&1& & & & & & \\
 & &3&3&3&1&1&1& & & & & & & & \\
 &3&3&1&1&1& & & & & & & & & & \\
 &1&1&1& & & & & & & & & & & & \\
1&1& & & & & & & & & & & & & &

\end{matrix}.
$$

\end{example}

We come now to the main result of this section.

\begin{theorem}
\label{reverse-ribbon-staircase-equivalence}
Let $\alpha$ be a  ribbon, and let
$$
\begin{aligned}
D &= (\varepsilon ^k _m(\alpha), \mathcal{N})_x\\
D' &= (\varepsilon ^k _m(\alpha), \mathcal{N}^\ast)_x
\end{aligned}
$$  
where $m < \ell(\alpha)$ and $x= se$ or $nw$.
Then $D \sim D'$.
\end{theorem}
\begin{proof} Assume that $x= se$; the case where $x=nw$ is analogous.

Index the ribbons in the southeast ribbon staircase decompositions of $D, D'$ 
so that the largest ribbon, which is the cutting strip $\theta$, comes first in each case.
Index the remaining ribbons so that they 
correspond under the natural bijection between the letters in the words
$\mathcal{N}$ and $\mathcal{N}^\ast$.  One can then check that the
associated Hamel-Goulden matrices are transposes of each other, and hence have the
same determinant.
\end{proof}

\begin{corollary}
\label{self-conjugate-corollary}
Let $D$ be a connected skew diagram with a ribbon staircase decomposition,
that is, $D=(\varepsilon ^k_m(\alpha), \mathcal{N})_x$ for some ribbon 
$\alpha$, with $m<l(\alpha)$ and $x= se$ or $nw$.  
Then $D^t$ also has a ribbon staircase decomposition, specifically
$$
D^t=(\varepsilon ^k _{m'}(\alpha^t), \mathcal{N}^*)_x
$$
where $m'= |\alpha \cap _m \alpha |-(m-1)$.  
Furthermore, if
$\alpha=\alpha^t$, then $m'=m$ and $D^t \sim D$.
\end{corollary}

\begin{proof}
The first assertion is a straightforward verification, in which one must
treat the cases $k=1,2$ separately.

For the second assertion, when $\alpha =\alpha ^t$ 
it is similarly straightforward to check that $m=m'$, and then one has
% Here is the long proof
%
%We claim
%$(\alpha \cap _m\alpha)^t$ = $\alpha ^t \cap _{m}\alpha ^t$ 
% for all $m$ for which $\alpha \cap _m \alpha$ is defined 
% and hence $m=m'$. If $\alpha = (\alpha _1,\ldots , \alpha _k)$ then there are three cases to consider.
%
%If $m=1$ then since $\alpha =\alpha ^t$ one of $\alpha _1,\alpha _k =1$ so
%$\alpha \cap _1 \alpha = \min \{\alpha _1, \alpha _k\}=1$ and
%$$(\alpha \cap _1 \alpha)^t= \alpha \cap _1 \alpha. $$
%
%If $m=2$ then observe since $\alpha =\alpha ^t$ then 
% one of $\alpha _1,\alpha _k =1$. If $\alpha _1 =1$ then $\alpha \cap _2\alpha $ 
% is undefined unless $\alpha _2=\alpha _k=2$. In this case $\alpha \cap _2\alpha =(1,2)$ and 
% $$(\alpha \cap _2 \alpha)^t= \alpha \cap _2 \alpha. $$Similarly,
% if $\alpha _k =1$ then $\alpha \cap _2\alpha $ is undefined 
% unless $\alpha _1=\alpha _{k-1}=2$. In this case $\alpha \cap _2\alpha =(2,1)$ and 
% $$(\alpha \cap _2 \alpha)^t= \alpha \cap _2 \alpha. $$
%
%Lastly, if $m>2$ then by definition the bottommost row  of $\alpha \cap _m \alpha$ 
% is the same length as the rightmost column of $\alpha \cap _m \alpha$, 
% the second from bottommost row  of $\alpha \cap _m \alpha$ is the same 
% length as the second from rightmost column of $\alpha \cap _m \alpha$, etc. 
% Hence it follows that $(\alpha \cap _m\alpha)^t=\alpha \cap _m \alpha$.
%
$$
D^t=(\varepsilon ^k _{m'}(\alpha^t), \mathcal{N}^*)_x
=(\varepsilon ^k _{m}(\alpha), \mathcal{N}^*)_x 
\sim D
$$
by Theorem~\ref{reverse-ribbon-staircase-equivalence}.  
\end{proof}

We close this section with an interesting special case of Corollary~\ref{self-conjugate-corollary}, which was first pointed out to us by John Stembridge and
for which we offer two proofs.

\begin{corollary}
\label{equivalent-conjugate-theorem}
For any Ferrers diagram $\mu$ contained in the staircase partition
$\delta_n:=(n-1,n-2,\ldots,1) \vdash \binom{n}{2},$
one has 
$$
\delta _n/\mu \sim \left( \delta _n / \mu \right) ^t.
$$
\end{corollary}

\vskip.1in
\noindent
{\it Proof 1.}
Check that the southeast decomposition of $\delta_n/\mu$ is always a 
southeast ribbon staircase decomposition of the form 
$\delta_n/\mu=(\varepsilon^{n-2}_1(\alpha),\mathcal{N})_{se}$,
in which $\alpha$ is the self-conjugate ribbon  $(1,2)$.  Then
apply Corollary~\ref{self-conjugate-corollary}.
\hfill$\qed$
% The detailed version of proof 1 is commented out here...
%
% We will prove that $\delta _n/\mu = (\varepsilon ^{n-2} _{1}((1,2)), \mathcal{N})_{se}$ 
% for some nesting $\mathcal{N}$. The result will then follow by 
% Corollary~\ref{self-conjugate-corollary}. We proceed by induction on $k\geq 3$.
%
%For $k=3$ we have $\delta _3=(\varepsilon ^1 _{1}((1,2)), \emptyset)_{se}$. 
%Now for $\delta _k/\mu$ where $\mu =(\mu _1,\ldots ,\mu _l)$, 
%consider the skew diagram $\delta _{k-1}/(\mu _2,\ldots ,\mu _l)$. 
%By induction we have  
%$$
%\delta _{k-1}/(\mu _2,\ldots ,\mu _l)=(\varepsilon ^{k-3} _{1}((1,2)), \mathcal{N})_{se}
%$$ 
%for some nesting $\mathcal{N}$. Let $k-2-\mu _1 = x$. Then depending on $x$ we now create a nesting 
%$\mathcal{N}'$ from $\mathcal{N}$ such that 
%$$
%\delta _k/(\mu _1,\ldots ,\mu _l) = (\varepsilon ^{k-2} _{1}((1,2)), \mathcal{N}')_{se}.
%$$
%
%Write $x=2q+r$ where $0\leq r <2$. If $\mu _1 >\mu _2$ then 
%if $r=0$ place a dot between the $q$ and $q+1$ letters of $\mathcal{N}$ 
%from the right, otherwise if $r=1$ then plcae a vertical slash between 
%the $q$ and $q+1$ letters of $\mathcal{N}$ from the right. If $\mu _1=\mu _2$ and 
%the rightmost dot or vertical slash in $\mathcal{N}$ is a dot then replace it 
%with $. |$, otherwise replace the vertical slash with $( )$ to create a new 
%nesting $\mathcal{N}'$. By construction it 
%follows $\delta _k/(\mu _1,\ldots ,\mu _l) = (\varepsilon ^{k-2} _{1}((1,2)), \mathcal{N}')_{se}$ and we are done.
\vskip.1in
\noindent
{\it Proof 2.} (cf. \cite[Proposition 7.17.7]{Stanley1})
Since all border strips in $D=\delta/\mu$ have odd size,
when one expands 
$$
s_{D}= \sum_{\lambda} z_\lambda^{-1}\chi^{D}(\lambda) p_\lambda
$$
in terms of power sum symmetric functions as in \cite[7.17.5]{Stanley1},
the Murnaghan-Nakayama formula \cite[Theorem 7.17.3]{Stanley1} for the coefficient
$\chi^{D}(\lambda)$ shows that it vanishes when $\lambda$ has any even parts.  Hence
$s_D$ is a polynomial in the odd power sums $p_1,p_3,p_5,\ldots$.
Since the involution $\omega$ satisfies $\omega(p_r)=(-1)^{r-1}p_r$, 
one has $s_{D^t}= \omega(s_D) = s_D$.\hfill
$\qed$

%%%%%%% Necessities start here %%%%%%%%

\section{Necessary conditions\label{necessities}}

 We now present some combinatorial invariants
for the skew-equivalence relation $D_1 \sim D_2$ on connected skew diagrams.
 
  \subsection{Frobenius rank}
  
The \emph{Durfee} or \emph{Frobenius rank} of a skew diagram $D$ is 
defined to be the minimum number of ribbons in any
decomposition of $D$ into ribbons.  The rank of $D$ is an invariant 
of $s_D$ which can be extracted in at least two ways using recent
results.  Firstly, Stanley has pointed out to us that
the discussion at the beginning of \cite[\S5]{Stanley2} implies
the rank of $D$ is the minimum length $\ell(\nu)$ among partitions
$\nu$ which appear when expanding $s_D$ uniquely as a sum of
power sum symmetric functions $p_\nu$.  Secondly,
it was recently conjectured by Stanley \cite{Stanley2}, and proven by
Chen and Yang \cite{ChenYang}, that the rank coincides
with the highest power of $t$ dividing the polynomial 
$s_D(1,1,\ldots,1,0,0,\ldots)$, where $t$ of the variables have been set 
to $1$, and the rest to zero.  Either of these implies the following.

\begin{corollary}
\label{zrank-corollary}
Frobenius rank is an invariant of skew-equivalence, that is  two skew-equivalent diagrams
must have the same Frobenius rank.

In particular, skew-equivalence restricts to the subset of
ribbons
as they are the skew diagrams of Frobenius rank $1$.
\end{corollary}

% Since having Durfee rank $1$ is the same as being a rim strip, which is the
% same as having overlap partition $\rho^{(2)}$ all ones, one might ask whether
% the Durfee rank can be recovered from the overlap partitions $\rho^{(k)}$;  I would
% assume not, but don't actually know ??

 \subsection{Overlaps}

Data about the amount of overlap between 
sets of rows or columns in the skew diagram $D$ can be recovered from
its skew Schur function $s_D$.

\begin{definition} \rm \ \\  
Let $D$ be a skew diagram occupying $r$ rows.  For each $k$ in $\{1,2,\ldots,r\}$,
define the \emph{$k$-row overlap composition} $r^{(k)}=(r^{(k)}_1,\ldots,r^{(k)}_{r-k+1})$ to be the sequence where
$r^{(k)}_i$ is the number of columns occupied in 
common by the rows $i,i+1,\cdots,i+k-1$.
Let $\rho^{(k)}$ be the \emph{$k$-row overlap partition} that is the weakly decreasing rearrangement
of $r^{(k)}$. Similarly define column overlap compositions $c^{(k)}$ and column overlap partitions $\gamma^{(k)}$.
\end{definition} 

\begin{example} \rm \ \\
If $D=\lambda/\mu$ with $\ell:=\ell(\lambda)$, then the $1$-row and $2$-row
overlap compositions are the sequences
$$
\begin{aligned}
r^{(1)} &= (\lambda_1-\mu_1,\ldots,\lambda_\ell-\mu_\ell) \\
r^{(2)} &= (\lambda_2-\mu_1,\lambda_3-\mu_2,\ldots,\lambda_{\ell}-\mu_{\ell -1})
\end{aligned}
$$
that played an important role in Proposition~\ref{JT-diagonals-proposition}.
\end{example}

It transpires that the row overlap partitions $( \rho^{(k)} )_{k \geq 1}$ and the
column overlap partitions $( \gamma^{(k)} )_{k \geq 1}$ determine each other uniquely.
To see this, we define a third form of data on a skew diagram $D$, which  mediates
between the two, and which is more symmetric under conjugation.

\begin{proposition}
\label{overlap-relation}Given a skew diagram $D$, consider the doubly-indexed array $( a_{k,\ell} )_{k,\ell \geq  1}$ where $a_{k,\ell}$ is defined
to be the number of $k \times \ell$ rectangular subdiagrams contained inside $D$.
Then we have
$$
\begin{aligned}
a_{k,\ell} &= \sum_{ \ell' \geq \ell } \left( \rho^{(k)} \right)^t_{\ell'} \\
           &= \sum_{ k' \geq k } \left( \gamma^{(\ell)} \right)^t_{k'}.
\end{aligned}
$$
Consequently, any one of the three forms of data
$$
( \rho^{(k)} )_{k \geq 1}, \quad
( \gamma^{(k)} )_{k \geq 1}, \quad
( a_{k,\ell} )_{k,\ell \geq  1}
$$
on $D$ determines the other two uniquely.

%In particular, one can recover $\gamma^{(1)}$ from $\rho^{(1)}, \rho^{(2)}, \ldots$ as follows:
%$$
%\left( \gamma^{(1)} \right)^t_k = | \rho^{(k)} | - | \rho^{(k-1)} |.
%$$
\end{proposition}
\begin{proof}
It suffices to prove the first equation, since exchanging rows and columns gives
the second. Every  $k \times \ell$ rectangular subdiagram of $D$ occupies
a particular $k$-tuple of rows, and the corresponding entry of $\rho^{(k)}$
coming from that $k$-tuple of rows must be of size $\ell' \geq \ell$.
This part $\ell'$ corresponds to a total of $\ell' - \ell + 1$ such
$k \times \ell$ subdiagrams, and hence
$$
\begin{aligned}
a_{k,\ell} &= \sum_{\text{ parts }\ell' \geq \ell\text{ in }\rho^{(k)}}
                 \left( \ell' -\ell + 1 \right) \\
           &= \sum_{ \ell' \geq \ell } \left( \rho^{(k)} \right)^t_{\ell'}.
\end{aligned}
$$
Since this relationship is invertible, the two forms of data
determine each other.
%For the last equation asserted in the proposition, note that both sides are
%counting (in different ways) the number of columns of $D$ of size at least $k$.
\end{proof}

Note that the data of the first two row-overlap compositions $r^{(1)}, r^{(2)}$ are enough to recover the
skew diagram $D$ up to translation within $\ZZ^2$, and similarly for
the compositions $c^{(1)}, c^{(2)}$.  
Thus one cannot expect to recover $r^{(1)}, r^{(2)}$ or $c^{(1)}, c^{(2)}$ from the
skew Schur function $s_D$.

However, it turns out 
one \emph{can} recover all of the row overlap {\it partitions} $( \rho^{(k)} )_{k \geq 1}$
(or the column overlap partitions $( \gamma^{(k)} )_{k \geq 1}$) from $s_D$;
see Corollary~\ref{overlaps-corollary} below, whose proof is our goal for the
remainder of this section.

\begin{lemma}
\label{hat-diagram-defined}
Given any skew diagram $D$, there is a unique skew diagram $\hat{D}$ satisfying
$r^{(k)}(\hat{D}) = r^{(k+1)}(D)$ for all $k$.
\end{lemma}
\begin{proof}
Observe that $\hat{D}$ is obtained from $D$ by removing the top cell from every column of $D$.
\end{proof}

\begin{example} \rm \ \\\label{hat-example}Let $D$ be the skew diagram from Example~\ref{LR-example}, shown below with its row overlap compositions 
depicted vertically to its right:
$$
D=\begin{matrix}
      &      &      &      & & r^{(1)}&r^{(2)}&r^{(3)}&r^{(4)}&r^{(5)}\\
      &      &      &\times& &  1&  &  &  &\\
      &      &      &      & &   &1 &  &  &\\
      &      &      &\times& &  1&  &0 &  &\\
      &      &      &      & &   &0 &  &0 &\\
      &\times&\times&      & &  2&  &0 &  &0\\
      &      &      &      & &   &2 &  &0 &\\
\times&\times&\times&      & &  3&  &1 &  &\\
      &      &      &      & &   &2 &  &  &\\
\times&\times&      &      & &  2&  &  &  &.
\end{matrix}
$$
Then $\hat{D}$ is as shown below: 
$$
\hat{D}=\begin{matrix}
      &      &      &      & &   &r^{(1)}&r^{(2)}&r^{(3)}&r^{(4)}\\
      &      &      &\times& &   &1 &  &  &\\
      &      &      &      & &   &  &0 &  &\\
      &      &      &      & &   &0 &  &0 &\\
      &      &      &      & &   &  &0 &  &0\\
      &\times&\times&      & &   &2 &  &0 &\\
      &      &      &      & &   &  &1 &  &\\
\times&\times&      &      & &   &2 &  &  &\\
\end{matrix}.
$$
Note in particular that the second row of $\hat{D}$ is {\it empty}.
\end{example}

We now prove a crucial lemma.

\begin{lemma}
\label{picture-lemma}
Let $D$ be a skew diagram with $c$ nonempty columns.
The map on column-strict tableaux that removes the first row and lowers all other entries by $1$ 
restricts to a bijection
$$
\{T \in \Pictures(D): \lambda_1(T)=c\} \longrightarrow \Pictures(\hat{D}).
$$
\end{lemma}
\begin{proof}
Let $T_{\row}(D)$ be the row filling of $D$, and let $t_j$ denote the topmost entry in
column $j$ of $T_{\row}(D)$ for $j=1,2,\ldots,c$.  
The following facts are then easy to verify from the definitions:
\begin{enumerate}

\item[(i)] the row filling $T_{\row}(\hat{D})$ is obtained from the row filling
$T_{\row}{(D)}$ by removing the entries $t_1,t_2,\ldots,t_c$ and then lowering the
remaining entries by $1$,
\item[(ii)] a picture $T$ for $D$ will have $\lambda_1(T)=c$ if and only if
the first row of $T$ is exactly $(t_c,t_{c-1},\ldots,t_2,t_1)$,
\item[(iii)] the map on pictures for $\hat{D}$ that raises all entries by $1$ 
and then adds a new top row $(t_c,t_{c-1},\ldots,t_2,t_1)$ gives a well-defined map
into the pictures for $D$, and is the inverse of the map defined in
the lemma.
\end{enumerate}
\end{proof}

\begin{example} \rm \ \\Let $D$ and $\hat{D}$ be as in Example~\ref{hat-example}, with row fillings
as shown here, and the topmost elements $(t_1,t_2,t_3,t_4)=(4,3,3,1)$ in each
column of $T_{\row}(D)$ shown in bold:
$$
T_{\row}(D)
=\begin{matrix}
 & & &\mathbf{1}\\
 & & &2\\
 &\mathbf{3}&\mathbf{3}& \\
\mathbf{4}&4&4& \\
5&5& & 
\end{matrix},
\qquad
T_{\row}(\hat{D})
=\begin{matrix}
 & & &1\\
 & & & \\
 &3&3& \\
4&4& &
\end{matrix}.
$$
Note that removing the entries $(t_1,t_2,t_3,t_4)$ in $T_{\row}(D)$ and lowering the
 remaining entries by $1$ gives $T_{\row}(\hat{D})$.
The set $\Pictures(D)$ was shown in Example~\ref{LR-example}, where
only the last five of these pictures $T$ have $\lambda_1(T)=c=4$:
$$
\left\{
\begin{matrix}
1&3&3&4\\
2&4& &\\
4&5& &\\
5& & &
\end{matrix},\quad 
\begin{matrix}
1&3&3&4\\
2&4&5&\\
4& & &\\
5& & &
\end{matrix},\quad 
\begin{matrix}
1&3&3&4\\
2&4&5&\\
4&5& &
\end{matrix},\quad 
\begin{matrix}
1&3&3&4\\
2&4&4&\\
5&5& &
\end{matrix},\quad 
\begin{matrix}
1&3&3&4\\
2&4&4&5\\
5& & &
\end{matrix} 
\right\}.
$$
On the other hand, for $\hat{D}$ one has pictures
$$
\left\{
\begin{matrix}
1&3& &\\
3&4& &\\
4& & &
\end{matrix},\quad 
\begin{matrix}
1&3&4&\\
3& & &\\
4& & &
\end{matrix},\quad 
\begin{matrix}
1&3&4&\\
3&4& &
\end{matrix},\quad 
\begin{matrix}
1&3&3&\\
4&4& &
\end{matrix},\quad 
\begin{matrix}
1&3&3&4\\
4& & &
\end{matrix} \right\}.
$$
\noindent
Note the bijection from the top set to the bottom set, obtained
by removing the first row and lowering the remaining entries by $1$.
\end{example}

\begin{definition} \rm \ \\
Given a positive integer $\ell$, define a $\ZZ$-linear map $\phi_\ell: \Lambda \rightarrow \Lambda$
by 
$$
\phi_\ell(s_\lambda) := 
\begin{cases}
   s_{\lambda+1^\ell} & \text{ if } \ell(\lambda) \leq \ell \\
   0                          & \text{ otherwise}
\end{cases}
$$
where 
$\lambda+1^\ell:=(\lambda_1+1,\ldots,\lambda_\ell+1)$.
Also, given $f \in \Lambda=\ZZ[h_1,h_2,\ldots]$, 
let $[h_r](f)$ denote the polynomial in the $h_i$ which gives the coefficient of $h_r$
in the expansion of $f$.
\end{definition}

\begin{theorem}
\label{hat-diagram-determined}
Let $D$ be a skew diagram with $\ell$ nonempty rows and $c$ nonempty columns.
Then 
$$
s_{\hat{D}} = [h_{\ell+c}] \phi_\ell(s_D).
$$
\end{theorem}
\begin{proof}
Recall the Littlewood-Richardson expansion \eqref{LR-rule}
$$
s_D = \sum_{T \in \Pictures(D)} s_{\lambda(T)}.
$$
Since $D$ has $\ell$ nonempty rows, any picture $T$ for $D$ will have at most
$\ell$ rows, and hence $\phi_\ell(s_{\lambda(T)})= s_{\lambda(T)+1^\ell}$.
Thus
$$
\begin{aligned}
\phi_\ell(s_D)
  &= \sum_{T \in \Pictures(D)} s_{\lambda(T)+1^\ell}\\
[h_{\ell+c}]\phi_\ell(s_D)
  &= \sum_{T \in \Pictures(D)} [h_{\ell+c}] s_{\lambda(T)+1^\ell}
\end{aligned}
$$
and it remains to extract the coefficient of $h_{\ell+c}$ in each term $s_{\lambda(T)+1^\ell}$.

It was noted in the proof of Corollary~\ref{reduction-to-connected-corollary} that
the Jacobi-Trudi expansion for $s_{\nu/\mu}$ takes a certain form;  when $\mu$ is empty
this form specializes to 
$$
s_\nu = s_{\hat{\nu}} \cdot h_{L(\nu)} + r
$$
where
$$
\begin{aligned}
L(\nu)   &:=\ell(\nu)+\nu_1-1\\
\hat{\nu}&:=(\nu_2-1,\nu_3-1,\ldots,\nu_{\ell(\nu)}-1)
\end{aligned}
$$
and the remainder $r$ is a polynomial in $h_k$ for $k < L(\nu)$.
Note that since $D$ has $c$ nonempty columns, any picture $T$ for $D$ will have
at most $c$ nonempty columns, and hence $\nu:=\lambda(T)+1^\ell$ will have
$L(\nu) \leq \ell+c$, with equality if and only if $\lambda_1(T)=c$.
Hence we have 
$$
[h_{\ell+c}] s_{\lambda(T)+1^\ell} =
\begin{cases}
s_{(\lambda_2(T),\ldots,\lambda_{\ell}(T))} & \text{ if } \lambda_1(T)=c \\
0                                             & \text{ otherwise.}
\end{cases}
$$
Consequently, 
$$
[h_{\ell+c}]\phi_\ell(s_D)
  = \sum_{\substack{T \in \Pictures(D):\\ \lambda_1(T)=c}} s_{(\lambda_2(T),\ldots,\lambda_{\ell}(T))},
$$
which equals $s_{\hat{D}}$ by Lemma~\ref{picture-lemma}.
\end{proof}

We are now ready to state our main necessary condition for skew-equivalence.
\begin{corollary}
\label{overlaps-corollary}
The skew Schur function $s_D$ determines the row overlap partition  $( \rho^{(k)} )_{k \geq 1}$ data.
Consequently, if $D \sim D'$, then $D, D'$ must have the same row overlap partitions.
\end{corollary}
\begin{proof}
Induct on the number $\ell$ of nonempty rows in $D$.  Proposition~\ref{JT-diagonals-proposition}(ii)
showed that $\rho^{(1)}$ can be recovered as the dominance-smallest partition occurring among the
subscripts of monomials in the $h_r$-expansion of $s_D$.

From Theorem~\ref{hat-diagram-determined} we know that $s_D$ determines $s_{\hat{D}}$.  By
induction, $s_{\hat{D}}$ determines its own row overlap partitions, which by Lemma~\ref{hat-diagram-defined}
coincide with the rest of the row overlap partitions $\rho^{(2)}, \rho^{(3)}, \ldots$ for $D$.
\end{proof}

\begin{example}  \rm \ \\Unfortunately, having the same row and column overlap partitions
$\rho^{(k)}, \gamma^{(k)}$ is not sufficient for the skew-equivalence of two skew diagrams.
For example, 
$$
\begin{matrix}
     &    &\times&\times \\
\times&\times&\times& \\
\times&      &      &
\end{matrix}
 \qquad \not\sim \qquad
\begin{matrix}
    &\times&\times&\times \\
\times&\times&      & \\
\times&      &      &
\end{matrix}
$$
even though they have the same row and column overlap partitions $\rho^{(k)}, \gamma^{(k)}$
for every $k$.
\end{example}

\begin{remark} \rm \ \\
Corollary~\ref{overlaps-corollary}
gives an alternate proof of the second assertion in Corollary~\ref{zrank-corollary},
that is, that ribbons can only be skew-equivalent to other ribbons, since
a skew diagram $D$ is a ribbon if and only if it is connected   and its $2$-row overlap partition $\rho^{(2)}$ has
the form $(1, 1, \ldots, 1)$.
\end{remark}

\section{Complete classification}\label{Conclusion}
 The sufficient conditions discussed in this paper explain all 
 but the following six skew-equivalences among skew diagrams with up to $18$ cells, up to antipodal rotation
and/or conjugation.  However, one way to explain both these skew-equivalences and the phenomenon occurring in Remark~\ref{smallest-nonribbon-circ-example} has recently been discovered in \cite{McNamaravW} and extends the definitions and results of Section~\ref{amalgamation-section} naturally.

$$
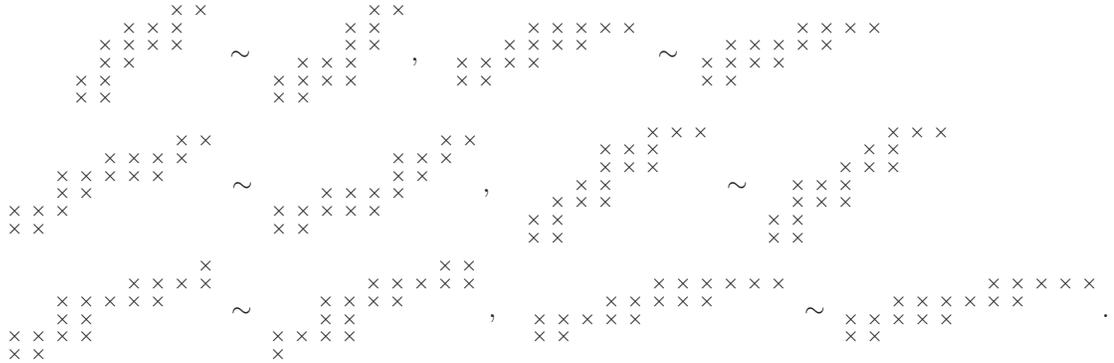
\begin{smallmatrix}
 & &  &    &\times&\times \\
 & &\times&\times&\times& \\
 &\times&     \times &   \times   &\times\\
 &\times&\times\\
\times&\times\\
\times&\times\end{smallmatrix}
 \  \sim \ 
\begin{smallmatrix}
 & &  &         &\times&\times \\
 & & &\times&\times& \\
 & & &   \times   &\times\\
 &\times&\times&\times\\
\times&\times&\times&\times\\
\times&\times\end{smallmatrix},\quad
%%%%%%%%%%%
 \begin{smallmatrix}
&&&\times&\times&\times&\times&\times\\
&&\times&\times&\times&\times\\
\times&\times&\times&\times\\
\times&\times\end{smallmatrix}
 \  \sim \ 
\begin{smallmatrix}
&&&&\times&\times&\times&\times\\
&\times&\times&\times&\times&\times\\
\times&\times&\times&\times\\
\times&\times\end{smallmatrix}$$
%%%%%%%%%%%%%%%
 $$\begin{smallmatrix}
&&&&&&&\times&\times\\
&&&&\times&\times&\times&\times&\\
&&\times&\times&\times&\times&\times\\
&&\times&\times\\
\times&\times&\times\\
\times&\times\end{smallmatrix}
 \  \sim \ 
\begin{smallmatrix}
&&&&&&&\times&\times\\
&&&&&\times&\times&\times\\
&&&&&\times&\times\\
&&\times&\times&\times&\times\\
\times&\times&\times&\times&\times\\
\times&\times\end{smallmatrix},\quad
%%%%%%%%%%%%
 \begin{smallmatrix}
&&&&&\times&\times&\times\\
&&&\times&\times&\times\\
&&&\times&\times&\times\\
&&\times&\times\\
&\times&\times&\times\\
\times&\times\\
\times&\times\end{smallmatrix}
 \ \sim \ 
\begin{smallmatrix}
&&&&&\times&\times&\times\\
&&&&\times&\times\\
&&&\times&\times&\times\\
&\times&\times&\times\\
&\times&\times&\times\\
\times&\times\\
\times&\times\end{smallmatrix}$$
%%%%%%%%%%%%%%%
$$\begin{smallmatrix}
&&&&&&&&\times\\
&&&&&\times&\times&\times&\times\\
&&\times&\times&\times&\times&\times\\
&&\times&\times\\
\times&\times&\times&\times&\\
\times&\times&\end{smallmatrix}
 \  \sim \ 
\begin{smallmatrix}
&&&&&&&\times&\times\\
&&&&\times&\times&\times&\times&\times&\\
&&\times&\times&\times&\times\\
&&\times&\times&\\
\times&\times&\times&\times&\\
\times&\end{smallmatrix},\quad
%%%%%%%%%%%%
\begin{smallmatrix}
 &&&&&\times&\times&\times&\times&\times&\times\\
 &&&\times&\times&\times&\times&\times\\
 \times&\times&\times&\times&\times\\
 \times&\times\end{smallmatrix} 
\  \sim \  
\begin{smallmatrix}
 &&&&&&\times&\times&\times&\times&\times\\
 &&\times&\times&\times&\times&\times&\times\\
 \times&\times&\times&\times&\times\\
 \times&\times\end{smallmatrix}.
$$

Note that these skew-equivalences occur in pairs. This leads us to end with the following conjecture that holds for all skew diagrams with up to 18 cells.

 \begin{conjecture} 
  Every skew-equivalence class of skew diagrams has cardinality a power of 2.
\end{conjecture}


\begin{thebibliography}{130}

\bibitem{BergeronBiagioliRosas}
F. Bergeron, R. Biagioli, and M.H. Rosas,
Inequalities between Littlewood-Richardson coefficients,
{\it J. Combin. Theory Ser. A} {\bf 113} (2006), 567-590.

\bibitem{BilleraThomasvanWilligenburg}
L.J. Billera, H. Thomas, and S. van Willigenburg,
Decomposable compositions, symmetric quasisymmetric functions and 
equality of ribbon Schur functions, {\it Adv.  Math.} {\bf204} (2006), 204--240.


\bibitem{Bressoud}
D. Bressoud,
Proofs and Confirmations: The Story of the Alternating-Sign Matrix Conjecture,
Cambridge University Press, Cambridge UK, 1999.

\bibitem{ChenYang}
W.Y.C. Chen and A.L.B. Yang,
Stanley's zrank problem on skew partitions,
to appear {\it Trans. Amer. Math. Soc.}.

\bibitem{ChenYanYang}
W.Y.C. Chen, G.-G. Yan, and A.L.B. Yang, 
Transformations of border strips and Schur function determinants,
{\it J. Algebraic Combin.} {\bf 21} (2005), 379--394.


\bibitem{FominFultonLiPoon}
S. Fomin, W. Fulton, C.-K. Li, and Y.-T. Poon,
Eigenvalues, singular values, and Littlewood-Richardson coefficients, 
{\it Amer. J.  Math.} {\bf 127} (2005), 101-127.

\bibitem{Gelfandetal}
I. M. Gel'fand, D. Krob, A. Lascoux, B. Leclerc, 
V. Retakh, and J.-Y. Thibon, Noncommutative symmetric functions, 
{\it Adv. Math.} {\bf 112} (1995), 218--348.

%\bibitem{FrumkinJamesRoichman}
%A. Frumkin, G. James, and Y. Roichman, 
%On trees and characters,
%{\it J. Algebraic Combin.} {\bf 17} (2003), 323--334.

\bibitem{HamelGoulden}
A. Hamel and I. Goulden,
Planar decompositions of tableaux and Schur function determinants,
{\it European J. Combin.} {\bf 16}  (1995), 461--477.

%\bibitem{Kraskiewicz}
%W. Kra\'skiewicz,
%Reduced decompositions in Weyl groups,
%{\it European J. Combin.} {\bf 16} (1995), 293--313.


\bibitem{LamPostnikovPylyavskyy}
T. Lam, A. Postnikov, and P. Pylyavskyy,
Schur positivity and Schur log-concavity,
to appear {\it Amer. J.  Math.}.

\bibitem{LascouxPragacz}
A. Lascoux and P. Pragacz,
\'Equerres et fonctions de Schur,
{\it C. R. Acad. Sci. Paris Sér. I Math.} {\bf 299} (1984), 955--958.

\bibitem{Macdonald} I.\ Macdonald, {\it Symmetric Functions and Hall Polynomials, 2nd Edition},
Oxford University Press, New York, USA, 1995.

%\bibitem{Magyar}
%P. Magyar,
%Borel-Weil theorem for configuration varieties and Schur modules,
%{\it Adv. Math.} {\bf 134} (1998), 328--366.

\bibitem{McNamaravW}
P.R.W. McNamara and S. van Willigenburg, A combinatorial classification of skew Schur functions, preprint {\tt math.CO/0608446}.

\bibitem{Okounkov}
A. Okounkov, Log-concavity of multiplicities with application to characters of $U(\infty)$,
{\it Adv. Math.} {\bf 127} (1997), 258-282.

%\bibitem{RhoadesSkandera}
%B. Rhoades and M. Skandera,
%Kazhdan-Lusztig immanants, to appear {\it J. Algebra}.

\bibitem{Sagan}
B.E. Sagan,
\emph{The symmetric group: Representations, combinatorial algorithms, and symmetric functions (2nd edition),
Graduate Texts in Mathematics, {\bf 203}}, Springer-Verlag, New York, USA, 2001. 

\bibitem{Schur}
I. Schur ``\"{U}ber eine Klasse von Matrizen die sich einer gegebenen Matrix zuorden lassen'', Inaugural-Dissertation, Berlin, 1901.

\bibitem{Shilov}
G.E. Shilov, 
\emph{Linear algebra. Revised English edition. Translated from the Russian and edited by Richard 
A. Silverman}, Dover Publications, New York, USA, 1977. 

\bibitem{Stanley1}
R.P. Stanley,
\emph{Enumerative Combinatorics, vol.~2},
Cambridge University Press, Cambridge, UK, 1999.

\bibitem{Stanley2}
R.P. Stanley,
The rank and minimal border strip decompositions of a skew partition, 
{\it J. Combin. Theory Ser. A} {\bf 100} (2002), 349--375.

%\bibitem{Stanley3}
%R.P.~Stanley,
%On the number of reduced decompositions of elements of Coxeter groups,
%{\it European J. Combin.} {\bf 5} (1984), 359--372.

\bibitem{YanYangZhou}
G.-G. Yan, A.L.B. Yang, and J. Zhou, 
The zrank conjecture and restricted Cauchy matrices,   {\it Linear Algebra Appl.}
{\bf 411} (2005), 371--385.

\bibitem{Zelevinsky}
A.V. Zelevinsky,
A generalization of the Littlewood-Richardson rule and the Robinson-Schensted-Knuth correspondence,
{\it J. Algebra} {\bf 69} (1981), 82--94.


\end{thebibliography}
\end{document}